\documentclass[11pt,twoside]{amsart}
\usepackage{amsmath, amsthm, amscd, amsfonts, amssymb, float, graphicx, color, xcolor, soul }
\usepackage[bookmarksnumbered, colorlinks, plainpages]{hyperref}
\usepackage{graphicx}
\usepackage{epstopdf}
\usepackage{rotating}
\input xy
\xyoption{all}
\sethlcolor{green}
\newtheorem{thm}{Theorem}[section]

\newtheorem{lem}[thm]{Lemma}
\newtheorem{prop}[thm]{Proposition}

\newtheorem{remark}{Remark}

\newtheorem{ex}[thm]{Example}

\numberwithin{equation}{section}

\def\Cay{{\mathrm {Cay}}}
\begin{document}
\openup 1.9\jot
% !TeX spellcheck = en_US
%------------------------------------------------------------------------------------%

\vspace{1.3 cm}
%------------------------------------------------------------------------------------%
\title{Domination parameters and diameter of Abelian Cayley graphs}
\author{Mohammad~A.~Iranmanesh$^*$ and Nasrin Moghaddami}

\thanks{{\scriptsize
\hskip -0.4 true cm MSC(2010): Primary: 20E99.
$*$Corresponding author}}
\maketitle
%------------------------------------------------------------------------------------%
%------------------------------------------------------------------------------------%
\begin{abstract}
Using the domination parameters of Cayley graphs constructed out of
$\mathbb{Z}_{p}\times \mathbb{Z}_{m}$, where $m\in\{p^{\alpha}, p^{\alpha}q^{\beta}, p^{\alpha}q^{\beta}r^{\gamma}\},$ in 
this paper we are discussing about the total and connected domination number and diameter of these Cayley graphs.

%\keywords{Cayley graph, total dominating set, connected dominating set, total domination number, connected domination number\\ \protect \indent 2010 {\it Mathematics Subject Classification.} Primary 20E99.}

\end{abstract}
\maketitle
\section{Introduction and Preliminaries}
\label{intro}
Let $(G, \cdot)$ be a group and $S=S^{-1}$ be a non empty subset of $G$ not containing the
identity element e of $G$. The simple graph $\Gamma$ whose vertex set $V(\Gamma)=G$ and edge set
$E(\Gamma)=\{\{v, vs\} | v\in V(\Gamma), s \in S \}$ is called the Cayley graph of
$G$ corresponding to the set $S$  and is denoted by $\Cay(G, S)$.
By $\mathbb{Z}_{n}$ we denote the abelian group of order $n$. For any vertex $v\in V(\Gamma)$,
the {\it open neighborhood} of $v$ is
the set $N(v)=\{u \in V (\Gamma) | \{u, v\}\in E(\Gamma)\}$ and the {\it closed neighborhood} of $v$
is the set $N[v]=N(v)\cup \{v\}$. For a set $X\subseteq V(\Gamma)$, the open neighborhood of
$X$ is $N(X)=\bigcup_{v\in X} N(v)$ and the closed neighborhood of $X$ is $N[X]=N(X)\cup X$ \cite{hay}.
A set $D\subseteq V(\Gamma)$ is said to be a {\it dominating set} if $N[D]=V(\Gamma)$ or equivalently,
every vertex in $V(\Gamma)\backslash D$ is adjacent to at least one vertex in $D$.
The {\it domination number} $\gamma(\Gamma)$ is the minimum cardinality of a
dominating set in $\Gamma$. A dominating set
with cardinality $\gamma(\Gamma)$ is called a {\it $\gamma$-set}.
 A set $T\subseteq V(\Gamma)$ is said to be a {\it total dominating set} if $N(T)=V(\Gamma)$ or equivalently,
every vertex in $V(\Gamma)$ is adjacent to a vertex in $T$.
The {\it total domination number} $\gamma_{t}(\Gamma)$ is the minimum cardinality of a total dominating 
set in $\Gamma$. A total dominating set
with cardinality $\gamma_{t}(\Gamma)$ is  called a {\it $\gamma_{t}$-set}.  A graph $\Gamma$ is said to be
 connected graph if there is at least one path between every pair of vertices in $\Gamma$. The connected components
  of a graph are its maximal connected subgraphs. A dominating set $D$ of $\Gamma$ is said to be a connected dominating
   set if the induced subgraph generated by $D$ is connected. The minimum cardinality of a connected dominating set of $\Gamma$
    is called the connected domination number of $\Gamma$ and is denoted by $\gamma_{c}(\Gamma)$, and the 
    corresponding set is denoted by {\it $\gamma_{c}$-set} of $\Gamma$. Let $\lambda$ be the length of the longest sequence
of consecutive integers in $\mathbb{Z}_{m}$, each of which shares a prime factor with $m$. Dominating sets were defined by Berge
and Ore \cite{berg,ore}. The concept of total domination in graphs was introduced by E.J. 
Cockayne and R.W. Dows and S.T. Hedetniemi \cite{daw}. S.T Hedetniemi, R.C. Laskar\cite{Hed} introduced the 
connected domination number in graphs. 
Madhavi \cite{mad} introduced the concept of Euler totient Cayley graphs and
their domination parameters studied by Uma Maheswary and B. Maheswary \cite{uma}.
Also some properties of direct product graphs of Cayley graphs with arithmetic graphs discussed by
Uma Maheswary and B. Maheswary \cite{uma2}, and their domination parameters studied by
Uma Maheswary and B. Maheswary and M. Manjuri \cite{uma1,uma3,yas}. 

A walk is a sequence of pairwise adjacent vertices of a graph. A path is a walk in which no vertex is repeated. 
The distance between two vertices of a graph is the number of edges of the shortest path between them. The
 diameter of a connected graph is the maximum distance between two vertices of the graph. Note that, the diameter of a 
 disconnected graph is considered, the maximum diameter of connected components of graph. Let $v, w \in V(\Gamma)$ then 
 the distance between $v, w$ is denoted by $d(v, w)$ and the diameter of $\Gamma$ is denoted by $diam(\Gamma)$ \cite{bi, bm}.
 
Here we study the total and connected dominating sets and diameter of Cayley
graphs constructed out of $\mathbb{Z}_{p}\times \mathbb{Z}_{m}$ where 
$m\in\{p^{\alpha}, p^{\alpha}q^{\beta}, p^{\alpha}q^{\beta}r^{\gamma}\}$.
The domination number of these graphs are presented in \cite{me} and we present some of the results without proofs .
\begin{thm}\label{p1}
Let $\Gamma=\Cay(\mathbb{Z}_{p}\times \mathbb{Z}_{p^{\alpha}}, \Phi)$. Then
\item[1)] $\gamma(\Gamma)=2$ where $p=2$ and $\alpha=1$.
\item[2)] $\gamma(\Gamma)=4$ where $p=2$ and $\alpha\geq 2$.
\item[3)] $\gamma(\Gamma)=3$ where $p\geq 3$ and $\alpha\geq 1$.
\end{thm}
\begin{thm}\label{p2}
Let $\Gamma=\Cay(\mathbb{Z}_{p}\times\mathbb{Z}_{p^{\alpha}q^{\beta}}, \Phi)$ where
 $p, q \geq 2$ and $\alpha, \beta\geq 1$. Then
 $\gamma(\Gamma)$ is given by Table~\ref{tab:01}.
 \begin{table}
\begin{center}
\caption{$\gamma(\Cay(\mathbb{Z}_{p}\times \mathbb{Z}_{p^{\alpha}q^{\beta}}, \Phi))$}\label{tab:01}
\begin{tabular}{|c|c|c|}
\hline
$\Gamma$& $\gamma(\Gamma)$ & Comments \\
\hline
$\Cay(\mathbb{Z}_{p}\times \mathbb{Z}_{pq}, \Phi)$ & $4$ &\\
\hline
$\Cay(\mathbb{Z}_{2}\times \mathbb{Z}_{2^{\alpha}q^{\beta}}, \Phi)$ & $8$ & $(\alpha, \beta)\neq (1, 1)$\\
\hline
$ \Cay(\mathbb{Z}_{p}\times \mathbb{Z}_{2^{\alpha}p^{\beta}}, \Phi)$ & $6$ & $(\alpha, \beta)\neq (1, 1)$\\
\hline
$ \Cay(\mathbb{Z}_{p}\times \mathbb{Z}_{p^{\alpha}q^{\beta}}, \Phi)$ & $5$ & $(\alpha, \beta)\neq (1, 1)$\\
&&$p=3, q\geq 5~or~q=3, p\geq 5$\\
\hline
$ \Cay(\mathbb{Z}_{p}\times \mathbb{Z}_{p^{\alpha}q^{\beta}}, \Phi)$ & $4$ & $(\alpha, \beta)\neq (1, 1)$\\
&&$p, q\geq 5$\\
\hline
\end{tabular}
\end{center}
\end{table}
\end{thm}

\begin{thm}\label{p3}
Let $\Gamma=\Cay(\mathbb{Z}_{p}\times\mathbb{Z}_{p^{\alpha}q^{\beta}r^{\gamma}}, \Phi)$ where
 $p, q, r \geq 2$ and $\alpha, \beta, \gamma \geq 1$. Then
 $\gamma(\Gamma)$ is given by Table~\ref{tab:02}.

\begin{table}
\begin{center}
\caption{$\gamma(\Cay(\mathbb{Z}_{p}\times\mathbb{Z}_{p^{\alpha}q^{\beta}r^{\gamma}}, \Phi))$}\label{tab:02}.
\begin{tabular}{|c|c|c|}
\hline
$\Gamma$& $\gamma(\Gamma)$ & Comments \\
\hline
$\Cay(\mathbb{Z}_{2}\times \mathbb{Z}_{2qr}, \Phi)$ & $8$ &\\
\hline
$\Cay(\mathbb{Z}_{p}\times \mathbb{Z}_{2pr}, \Phi)$ & $8$ &  \\
\hline
$ \Cay(\mathbb{Z}_{2}\times \mathbb{Z}_{2^{\alpha}q^{\beta}r^{\gamma}}, \Phi)$ & $12$ & $\alpha\neq 1~or~\beta\neq 1~or~\gamma\neq 1$\\
\hline
$ \Cay(\mathbb{Z}_{p}\times \mathbb{Z}_{2^{\alpha}p^{\beta}r^{\gamma}}, \Phi)$ & $10$ & $\alpha\neq 1~or~\beta\neq 1~or~\gamma\neq 1$\\
&&$p=3, r\geq 5~or~r=3, p\geq 5$\\
\hline
$ \Cay(\mathbb{Z}_{p}\times \mathbb{Z}_{2^{\alpha}p^{\beta}r^{\gamma}}, \Phi)$ & $8$ & $\alpha\neq 1~or~\beta\neq 1~or~\gamma\neq 1$\\
&&$p,r\geq 5$\\
\hline
$ \Cay(\mathbb{Z}_{p}\times \mathbb{Z}_{p^{\alpha}q^{\beta}r^{\gamma}}, \Phi)$ & $6\leq\gamma(\Gamma)\leq 8$ & $\alpha, \beta, \gamma\geq 1$ \\
&& one of the prime factors is $3$\\
\hline
$\Cay(\mathbb{Z}_{p}\times \mathbb{Z}_{p^{\alpha}q^{\beta}r^{\gamma}}, \Phi)$ & $5$ & $p, q, r \geq 5 $ and $\alpha, \beta, \gamma\geq 1$\\
\hline
\end{tabular}
\end{center}
\end{table}
\end{thm}
Let $ p_{1}, p_{2}, \ldots, p_{k}$ be consecutive prime numbers, $\alpha, \alpha_{1}, \alpha_{2}, \ldots, \alpha_{k}$ are
positive integers and $\Phi=\varphi_{2}\times \varphi_{2^{\alpha}p_{1}^{\alpha_{1}}p_{2}^{\alpha_{2}}
 \ldots p_{k}^{\alpha_{k}}}$.

\begin{thm}\label{p4}
Let $\Gamma=\Cay(\mathbb{Z}_{2}\times \mathbb{Z}_{2^{\alpha}p_{1}^{\alpha_{1}}p_{2}^{\alpha_{2}}
 \ldots p_{k}^{\alpha_{k}}}, \Phi)$, where $p_{1}=3$ and $\alpha\geq 2$.
  Then $\gamma(\Gamma)\geq 4k+4$.
 \end{thm}

For $p=2$, the Cayley graph $\Cay(\mathbb{Z}_{p}\times \mathbb{Z}_{m},\Phi)$,
where $\Phi=\varphi_{p}\times \varphi_{m}$, is a disconnected graph with two connected
components, say $\Gamma_{1}$ and $\Gamma_{2}$, where
$V(\Gamma_{1})=\{(1, v)| v~is~odd\}\cup \{(0, v)| v~is~even\}$ and
$V(\Gamma_{2})=\{(0, v)| v~is~odd\}\cup \{(1, v)| v~is~even\}$.
Since every Cayley graph $\Cay(G, S)$ is $|S|$-regular (see for example \cite{godsil}),
we find that $\Gamma$ is $|\Phi|$-regular.

 Let $X$ be a set of 
consecutive integers in $\mathbb{Z}_{m}$ such that for
every $x\in X$, we have $\gcd(x, m)>1$. In this case we call $X_{i}$ a consecutive set. We use $X_{i}^{k}$ to show 
that the consecutive set $X_{i}$ has $k$ elemens.

Let $\Gamma=\Cay(\mathbb{Z}_{p}\times \mathbb{Z}_{m},\Phi)$. In Section \ref{sec:2}
we calculate $\gamma_{t}(\Gamma)$ and $\gamma_{c}(\Gamma)$ and $diam(\Gamma)$ where $m=p^{\alpha}$. 
We consider the case $m=p^{\alpha}q^{\beta}$ in Section \ref{sec:3} and the case
$m=p^{\alpha}q^{\beta}r^{\gamma}$ is considered in Section \ref{sec:4}.

\section{Total and connected domination number and diameter of $\Cay(\mathbb{Z}_{p}\times \mathbb{Z}_{p^{\alpha}}, \Phi)$}\label{sec:2}
Let $p$ be a prime number, $\alpha$ a positive integer and $\Phi=\varphi_{p}\times \varphi_{p^{\alpha}}$.
In this section, we obtain the total and connected domination number and diameter of
$\Gamma=\Cay(\mathbb{Z}_{p}\times \mathbb{Z}_{p^{\alpha}}, \Phi)$.

\begin{thm}\label{d1}
Let $\Gamma=\Cay(\mathbb{Z}_{p}\times \mathbb{Z}_{p^{\alpha}}, \Phi)$. Then
\item[1)]$diam(\Gamma)=1$ where $p=2$ and $\alpha=1$.
\item[2)]$diam(\Gamma)=2$ where $p=2, \alpha\geq 2$ and $p\geq 3, \alpha\geq 1$
\begin{proof}
$1)$ In this case $\Gamma\cong 2K_{2}$, and clearly the diameter of $\Gamma$ is $1$.

$2)$ Let $p=2$ and $\alpha\geq 2$. Then $\Gamma$ is a disconnected graph with two connected
components, say $\Gamma_{1}$ and $\Gamma_{2}$, where
$V(\Gamma_{1})=\{(1, v)| v~is~odd\}\cup \{(0, v)| v~is~even\}$ and
$V(\Gamma_{2})=\{(0, v)| v~is~odd\}\cup \{(1, v)| v~is~even\}$.

Let $(u, v), (u^{'}, v^{'})\in V(\Gamma_{1})$. Then we have the following two possibilities:

$i)$ $u=u^{'}$ and $v\neq v^{'}$. Obviously $(u, v)$ and $(u, v^{'})$ are not adjacent. This implies that 
$d((u, v), (u^{'}, v^{'}))\geq 2$. On the other hand the vertex $(u-1, v-1)$ is adjacent to both vertices. So $d((u, v), (u^{'}, v^{'}))=2$. 

$ii)$ $u\neq u^{'}$ and $v\neq v^{'}$. We know that $u- u^{'} \in \varphi_{2}$ and $v-v^{'}$ is an odd integer.
 Since all of the odd integers in $\mathbb{Z}_{2^{\alpha}}$ to be included into a $\varphi_{2^{\alpha}}$, hence
  $v-v^{'}\in \varphi_{2^{\alpha}}$. Thus $(u, v)$ is adjacent to $(u^{'}, v^{'})$. So  $d((u, v), (u^{'}, v^{'}))=1$.

Since $(u, v)$ and $(u^{'}, v^{'})$ are arbitrary vertices of $\Gamma_{1}$, hence the diameter of $\Gamma_{1}$ is $2$.
 Similarly the diameter of $\Gamma_{2}$ is $2$. 
 
 Let $p\geq 3$ and $\alpha\geq 1$. Then $\Gamma$ is connected graph where 
 $$V(\Gamma)=\{(0, 0), \ldots, (0, p^{\alpha}-1), \ldots, (p-1,0), \ldots, (p-1,p^{\alpha}-1)\}.$$
 Assume that $(u, v)$ and $(u^{'}, v^{'})$ are arbitrary vertices of $\Gamma$.
  Now we have the following three possibilities:
 
 $i)$ $u=u^{'}$ and $v\neq v^{'}$. Since $(u, v)$ and $(u^{'}, v^{'})$ are not adjacent, Hence $d((u, v), (u^{'}, v^{'}))\geq 2$. 
  Let $v$ and $v^{'}$ be multiple of $p$. Note that $0$ is multiple of $p$. Then $(u-1, p-1)$ is adjacent to both 
 $(u, v)$ and $(u^{'}, v^{'})$. Let $v$ and $v^{'}$ be non-multiple of $p$. Then $(u-1, p)$ is common neigbor of 
 $(u, v)$ and $(u^{'}, v^{'})$. Now let one of either $v$ or $v^{'}$ is multiple of $p$. Without loss of generality let $v$ is multiple of $p$ and $v^{'}$ is non-multiple of $p$. Suppose that $v$ and $v^{'}$ are both even or odd.
 Then $(u-1, \dfrac{v+v^{'}}{2})$ is adjacent to both $(u, v)$ 
 and $(u^{'}, v^{'})$. Since $v-v^{'}$ is even so $v-v^{'}$ is 
 divisible by $2$. Hence  $v-\dfrac{v+v^{'}}{2}=\dfrac{2v-v-v^{'}}{2}=\dfrac{v-v^{'}}{2}\in\varphi_{p^{\alpha}}$ and also 
 $v^{'}-\dfrac{v+v^{'}}{2}=\dfrac{v^{'}-v}{2}\in\varphi_{p^{\alpha}}$. Now assume that one of either 
 $v$ or $v^{'}$ is even. Then $(u-1, 2v^{'})$ is
   common neigbor of $(u, v)$ and $(u^{'}, v^{'})$. Therefore $d((u, v), (u^{'}, v^{'}))=2$.
  
  $ii)$ $u\neq u^{'}$ and $v=v^{'}$. In this case vertex $(u^{''}, v-1)$ where $u^{''}\neq u, u^{'}$ is adjacent 
  to both $(u, v)$ and $(u^{'}, v^{'})$. Thus $d((u, v), (u^{'}, v^{'}))=2$.
  
  $iii)$ $u\neq u^{'}$ and $v\neq v^{'}$. If $(u, v)$ and $(u^{'}, v^{'})$ be adjacent then $d((u, v), (u^{'}, v^{'}))=1$. 
  If $(u, v)$ and $(u^{'}, v^{'})$ be non-adjacent then by $(i)$ and $(ii)$, $d((u, v), (u^{'}, v^{'}))=2$. Therefore in this case $diam(\Gamma)=2$.
\end{proof}
\end{thm}
\begin{thm}
Let $\Gamma=\Cay(\mathbb{Z}_{p}\times \mathbb{Z}_{p^{\alpha}}, \Phi)$. Then
\item[1)]$\gamma_{t}(\Gamma)=4$ and $\gamma_{c}(\Gamma)$ does not exist where $p=2$ and $\alpha\geq 1$.
\item[2)]$\gamma_{t}(\Gamma)=\gamma_{c}(\Gamma)=3$ where $p\geq 3$ and $\alpha\geq 1$.
\begin{proof}
 $1)$ Let $p=2$ and $\alpha=1$. Then $\Gamma\cong 2K_{2}$, and obviously
$\gamma_{t}(\Gamma)=4$.

Assume that  $p=2$ and $\alpha\geq 2$. Then by \cite[Theorem 2.1]{me}, $\gamma(\Gamma)=4$ and 
$D=\{(0, 0), (0, 1), \\(1, 0), (1, 1)\}$ is a $\gamma$-set for $\Gamma$. Since $(0, 0)$ and $(0, 1)$ are 
adjacent to $(1, 1)$ and $(1, 0)$, respectively hence $D$ is a 
$\gamma_{t}$-set for $\Gamma$. Thus $\gamma_{t}(\Gamma)=4$.

 In this case $\Gamma$ is a disconnected graph. 
Hence by the definition of connected dominating set, {\it $\gamma_{c}$-set} does not exist for $\Gamma$

$2)$ Let $p\geq 3$ and $\alpha\geq 1$. By \cite[Theorem 2.1]{me}, $\gamma(\Gamma)=3$ 
and $D=\{(0, 1), (1, 0),(2, 2)\}$ is a $\gamma$-set for $\Gamma$. Vertices of $D$ dominate 
among themselves. Therefore $\gamma_{t}(\Gamma)=\gamma_{c}(\Gamma)=3$. 
\end{proof}
\end{thm}

\begin{figure}[h]
\begin{center}
\includegraphics[width=14.5cm , height=6.5cm]{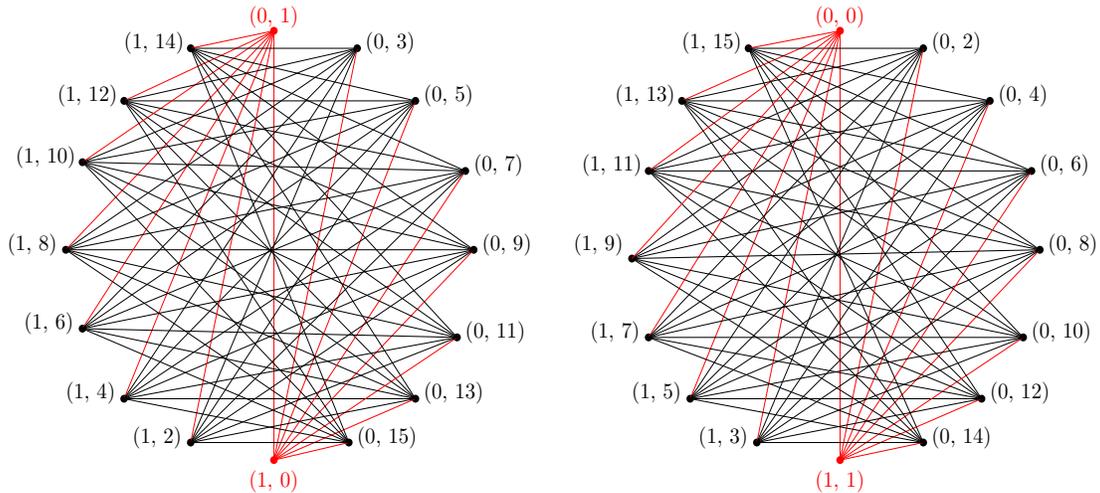}
%\vspace*{-.2cm}
\caption{The graph $\Gamma_{1}=\Cay(\mathbb{Z}_{2}\times \mathbb{Z}_{2^{4}},\Phi)$ and its total dominating set.}\label{fig:1}
\end{center}
\end{figure}

\begin{figure}[h]
\begin{center}
\includegraphics[width=8.5cm , height=6cm]{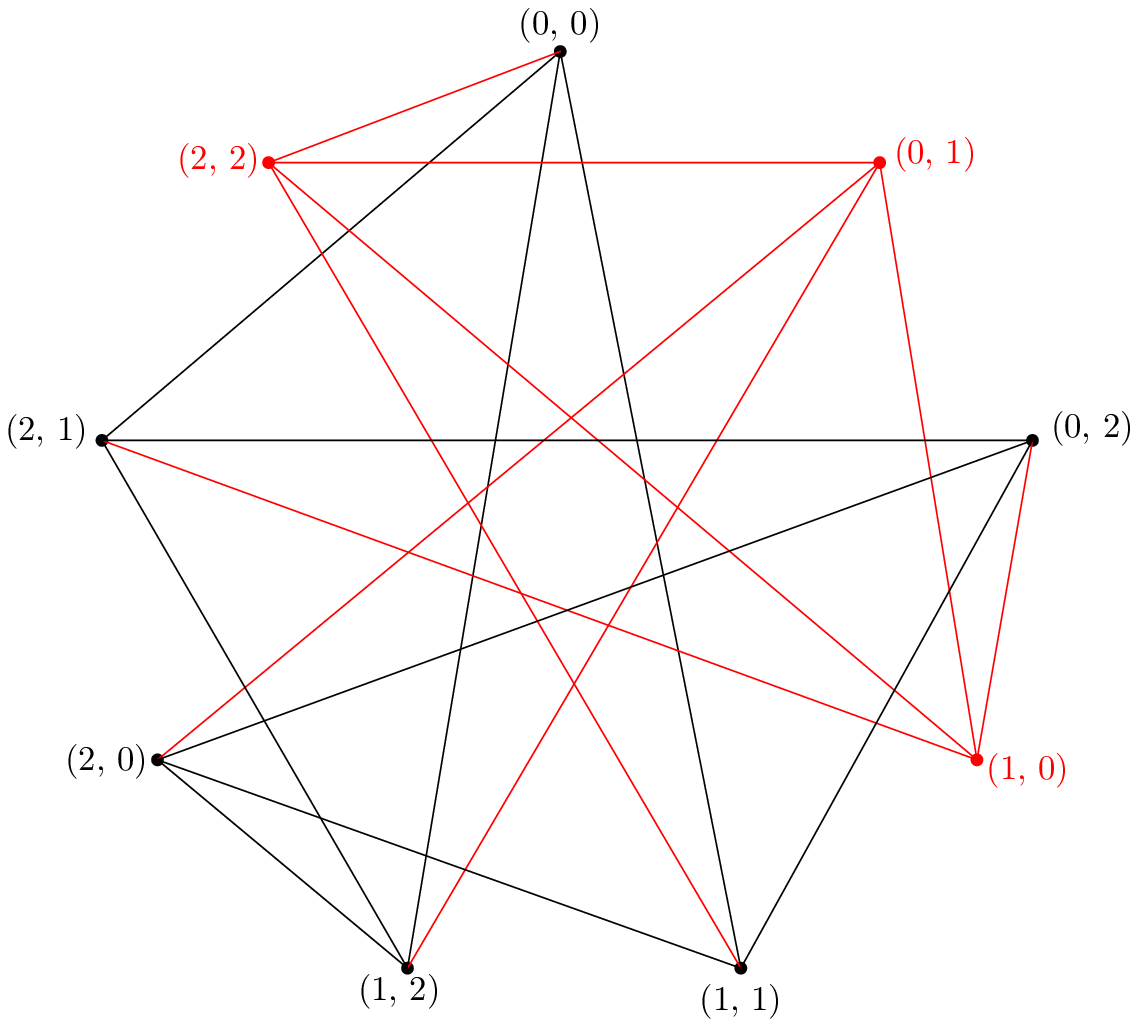}
%\vspace*{-.2cm}
\caption{The graph $\Gamma_{2}=\Cay(\mathbb{Z}_{3}\times \mathbb{Z}_{3},\Phi)$ and
 its total and connected dominating set.}\label{fig:2}
\end{center}
\end{figure}

\begin{ex}
Let $\Gamma_{1}=\Cay(\mathbb{Z}_{2}\times \mathbb{Z}_{2^{4}},\Phi)$
and $\Gamma_{2}=\Cay(\mathbb{Z}_{3}\times \mathbb{Z}_{3},\Phi),$
which are shown in Figures~\ref{fig:1} and~\ref{fig:2}, respectively. Clearly $\Gamma_{1}$
is a disconnected graph with two connected components. Thus {\it $\gamma_{c}$-set} does not 
exist for $\Gamma_{1}$. Also total dominating set of $\Gamma_{1}$, is $\{(0, 0), (0, 1), (1, 0), (1, 1)\}$. Note that
 total and connected dominating set of $\Gamma_{2}$ is $\{(0, 1), (1, 0), (2, 2)\}$.
\end{ex}


\section{Total and connected domination number and diameter of $\Cay(\mathbb{Z}_{p}\times \mathbb{Z}_{p^{\alpha}q^{\beta}},\Phi)$}\label{sec:3}
Let $p, q$ be prime numbers, $\alpha, \beta$ positive integers and $\Phi=\varphi_{p}\times \varphi_{p^{\alpha}q^{\beta}}$.
In this section, we find the total and connected domination number and diameter of
$\Gamma=\Cay(\mathbb{Z}_{p}\times \mathbb{Z}_{p^{\alpha}q^{\beta}},\Phi)$.
\begin{lem}\label{d2}
Let $\Gamma=\Cay(\mathbb{Z}_{2}\times \mathbb{Z}_{2^{\alpha}q^{\beta}},\Phi)$, where $\alpha, \beta \geq 1$. Then $diam(\Gamma)=3$.
\begin{proof}
$\Gamma$ is a disconnected graph with two connected
components, say $\Gamma_{1}$ and $\Gamma_{2}$, where
$V(\Gamma_{1})=\{(1, v)| v~is~odd\}\cup \{(0, v)| v~is~even\}$ and
$V(\Gamma_{2})=\{(0, v)| v~is~odd\}\cup \{(1, v)| v~is~even\}$.

Let $(u, v), (u^{'}, v^{'})\in V(\Gamma_{1})$. Then we have the following two possibilities:

$i)$ $u=u^{'}$ and $v\neq v^{'}$. Clearly $d((u, v), (u^{'}, v^{'}))\geq 2$. Let $v$ and $v^{'}$ be
 multiple of $2q$. Then $(u-1, 2q-1)$ is common neighbor of $(u, v)$ and  $(u^{'}, v^{'})$. Also if $v$ 
 and $v^{'}$ be non-multiple of $2q$, then $(u-1, 2q)$ is adjacent to both $(u, v)$ and  $(u^{'}, v^{'})$. Let $v$ and $v^{'}$ be both 
 multiple of one of the prime factors $2$ or $q$. Then the other prime factor is adjacent to both $v$ 
 and $v^{'}$. Now let one of either $v$ or $v^{'}$ is odd and is multiple of $q$. Then $(u, v), (u^{'}, v^{'})\in \{(1, v)| v~is~odd\}$. If $\dfrac{v+v^{'}}{2}$ be even, then
  $(u-1, \dfrac{v+v^{'}}{2})$ is common neighbor of $(u, v)$ and  $(u^{'}, v^{'})$. Also if $\dfrac{v+v^{'}}{2}$ be odd, then
  $(u-1, \dfrac{v+v^{'}}{2}+q)$ is adjacent to both $(u, v)$ and  $(u^{'}, v^{'})$. Let one of either $v$ or $v^{'}$
   is multiple of $2q$. So $(u, v), (u^{'}, v^{'})\in \{(0, v)| v~is~even\} $. If $\dfrac{v+v^{'}}{2}\in \varphi_{2^{\alpha}q^{\beta}}$, then
  $(u-1, \dfrac{v+v^{'}}{2})$ is common neighbor of $(u, v)$ and  $(u^{'}, v^{'})$. More over if $\dfrac{v+v^{'}}{2}$ be even, then
  $(u-1, \dfrac{v+v^{'}}{2}+q)$ is adjacent to both $(u, v)$ and  $(u^{'}, v^{'})$. Thus in this case $d((u, v), (u^{'}, v^{'}))=2$.
  
 $ii)$ $u\neq u^{'}$ and $v\neq v^{'}$. If $v-v^{'}\in \varphi_{2^{\alpha}q^{\beta}}$, then $d((u, v), (u^{'}, v^{'}))=1$. 
  Suppose that $v-v^{'}\notin \varphi_{2^{\alpha}q^{\beta}}$, since $u\neq u^{'}$ and $u, u^{'}\in \mathbb{Z}_{2}$, 
  we have no common neighbor between $(u, v)$ and  $(u^{'}, v^{'})$. This implies that $d((u, v), (u^{'}, v^{'}))\geq 3$. 
   Without loss of generality assume that $u=0$ and $u^{'}=1$. Since $v-v^{'}$ is an odd integer, we find that $v-v^{'}-2\in \varphi_{2^{\alpha}q^{\beta}}$. Thus $(0, v)(1, v+1)(0, v+2)(1, v^{'})$ is a 
   path of length $3$ between $(0, v)$ and  $(1, v^{'})$. So $diam(\Gamma_{1})=3$ and similarly $diam(\Gamma_{2})=3$. Therefore $diam(\Gamma)=3$.
\end{proof}
\end{lem}

\begin{lem}\label{t1}
Let $\Gamma=\Cay(\mathbb{Z}_{2}\times \mathbb{Z}_{2^{\alpha}q^{\beta}},\Phi)$, where $\alpha, \beta \geq 1$. 
Then $\gamma_{c}(\Gamma)$ does not exist and $\gamma_{t}(\Gamma)=8$.

\begin{proof}
$\Gamma$ is a disconnected graph with exactly two
connected components $\Gamma_{1}$ and $\Gamma_{2}$ where $V(\Gamma_{1})=\{(1, v)| v~is~odd\}\cup \{(0, v)| v~is~even\}$ and
$V(\Gamma_{2})=\{(0, v)| v~is~odd\}\cup \{(1, v)| v~is~even\}.$ Hence by the definition of 
connected dominating set, {\it $\gamma_{c}$-set} does not exist for $\Gamma$. 

Assume first that $(\alpha, \beta)=(1, 1)$. Then by \cite[Proposition 3.1]{me}, $A=\{(0, 0), (1, q)\}$ 
and $B=\{(0, 1), (1, q+1)\}$ dominate $V(\Gamma_{1})\backslash A$ and $V(\Gamma_{2})\backslash B$, respectively.
 Hence $\gamma(\Gamma)=4$. Vertices
  of $A$ are not adjacent to each other and  $A$ is not dominated by one vertex. Note that  $(1, 1)$ and 
  $(0, q+1)$ are adjacent to $(0, 0)$ and $(1, q)$, respectively. Hence
   $T_{1}=\{(0, 0), (1, 1), (1, q), (0, q+1)\}$ is a $\gamma_{t}$-set for $\Gamma_{1}$. 
   Similarly $T_{2}=\{(0, 1), (1, 0), (0, q), (1, q+1)\}$ is a $\gamma_{t}$-set for $\Gamma_{2}$. Therefore $\gamma_{t}(\Gamma)=8$.
 
Next consider the case where $(\alpha, \beta)\neq (1, 1)$. By \cite[Lemma 3.2]{me}, $\gamma(\Gamma)=8$ 
and $D=\{(0, 0), (0, 1), (0, 2), (0, 3), (1, 0), (1, 1), (1, 2), (1, 3)\}$ is a $\gamma$-set for $\Gamma$. 
Vertices $(0, 1), (0, 0), (0, 3),\\ (0, 2)$ are adjacent to vertices $(1, 0), (1, 1), (1, 2), (1, 3)$ respectively. Thus $D$ becomes
 a $\gamma_{t}$-set for $\Gamma$. Hence $\gamma_{t}(\Gamma)=8$.
 \end{proof}
\end{lem}
\begin{prop}\label{d3}
Let $\Gamma=\Cay(\mathbb{Z}_{p}\times \mathbb{Z}_{2^{\alpha}p^{\beta}},\Phi)$, 
where $\alpha, \beta \geq 1$. Then $diam(\Gamma)=3$.
\begin{proof}
Let $(u, v), (u^{'}, v^{'})\in V(\Gamma)$. Then we have the following three possibilities:

$i)$ $u=u^{'}$ and $v\neq v^{'}$. In this case $d((u, v), (u^{'}, v^{'}))\geq 2$. Suppose that $v$ and $v^{'}$ 
are both even or odd. Hence by case $(i)$ of Theorem~\ref{d2}, $d((u, v), (u^{'}, v^{'}))=2$. Since in 
$\mathbb{Z}_{2}\times \mathbb{Z}_{2^{\alpha}p^{\beta}}$ we have two connected components, where in 
each of them, if $u=u^{'}$ then $v$ and $v^{'}$ are both even or odd. 

Assume that one of either $v$ or $v^{'}$
 is even. Without loss of generality let $v$ is even and $v^{'}$ is odd. Also 
 let $(u^{''}, v^{''})$ where $u^{''}\neq u$, is common neighbor between $(u, v), (u^{'}, v^{'})$. If 
  $v^{''}$ be even then $v-v^{''}\notin\varphi_{2^{\alpha}p^{\beta}}$ and if $v^{''}$ be odd then
   $v^{'}-v^{''}\notin\varphi_{2^{\alpha}p^{\beta}}$. Thus we have no common neighbor between
    $(u, v)$ and $(u^{'}, v^{'})$. Hence $d((u, v), (u^{'}, v^{'}))\geq 3$. We consider $u^{''}, u^{'''}\neq u$, if $v$ and $v^{'}$
   be multiple of $p$, then $(u, v) (u^{''}, p-2) (u^{'''}, p-1) (u, v^{'})$ is a path 
  of length $3$ between $(u, v)$ and $(u^{'}, v^{'})$. If  $v$ and $v^{'}$
   be non-multiple of $p$, then the path $(u, v) (u^{''}, p) (u^{'''}, 2v^{'}) (u^{'}, v^{'})$ is connected.
   If $v$ be multiple of $p$ and $v^{'}$
   be non-multiple of $p$, since $v-v^{'}\in \varphi_{2^{\alpha}p^{\beta}}$ then $(u, v) (u^{''}, v^{'}) (u^{'''}, v) (u^{'}, v^{'})$
   is a path of length 3 between $(u, v)$ and $(u^{'}, v^{'})$.
   
 $ii)$ $u\neq u^{'}$ and $v=v^{'}$. In this case $(u^{''}, v-1)$ where $u^{''}\neq u, u^{'}$ is common neighbor of $(u, v)$ 
 and $(u^{'}, v^{'})$. Hence $d((u, v), (u^{'}, v^{'}))=2$.
 
 $iii)$ $u\neq u^{'}$ and $v\neq v^{'}$. If $(u, v)$ and $(u^{'}, v^{'})$ be adjacent then $d((u, v), (u^{'}, v^{'}))=1$. 
 Now assume that  $(u, v)$ and $(u^{'}, v^{'})$ are not adjacent. Let  $v$ and $v^{'}$ 
be both even or odd. Then by case $(i)$ of Lemma~\ref{d2}, we know that there is a vertex $(u^{''}, v^{''})$, 
where $u^{''}\neq u, u^{'}$ and $v^{''}$ is adjacent to $v$ and $v^{'}$, that is adjacent to both $(u, v)$ and $(u^{'}, v^{'})$. Thus $d((u, v), (u^{'}, v^{'}))=2$. 

Now let one of either $v$ or $v^{'}$ is even. Then by second paragraph of case $(i)$ and also by using of case $(ii)$ of 
Lemma~\ref{d2}, we see that $d((u, v), (u^{'}, v^{'}))=3$. Therefore $diam(\Gamma)=3$.
\end{proof}
\end{prop}
\begin{prop}\label{t2}
Let $\Gamma=\Cay(\mathbb{Z}_{p}\times \mathbb{Z}_{2^{\alpha}p^{\beta}},\Phi)$, where $\alpha, \beta \geq 1$. 
Then
\item[1)] $\gamma_{t}(\Gamma)=6$.
\item[2)] $\gamma_{c}(\Gamma)$ is given by Table~\ref{tab:7}.
\begin{table}
\begin{center}
\caption{$\gamma_{c}(\Cay(\mathbb{Z}_{p}\times \mathbb{Z}_{2^{\alpha}p^{\beta}}, \Phi))$ where $\alpha, \beta\geq 1.$}\label{tab:7}
\begin{tabular}{|c|c|c|}
\hline
$\Gamma$& $\gamma_{c}(\Gamma)$ & Comments \\
\hline
$\Cay(\mathbb{Z}_{p}\times \mathbb{Z}_{2^{\alpha}p^{\beta}}, \Phi)$ & $7$ & $p=3$\\
\hline
$\Cay(\mathbb{Z}_{p}\times \mathbb{Z}_{2^{\alpha}p^{\beta}}, \Phi)$ & $6$ & $p \geq 5 $  \\
\hline
\end{tabular}
\end{center}
\end{table}
\begin{proof}
 $1)$ Let $(\alpha, \beta)=(1, 1)$. By \cite[Proposition 3.1]{me}, $\gamma(\Gamma)=4$ and $D=\{(0,0), (0, 1), (1, p), (1,\\ p+1)\}$
 is a $\gamma$-set for $\Gamma$. Vertices
  of $D$ are not adjacent to each other. Hence $\gamma_{t}(\Gamma)>4$. Let a vertex say $(u, v)$ dominates all vertices of $D$. 
  Then $(u, v)$ is adjacent to $(0, 0)$ hence $(u, v)\in \Phi$. On the other hand $(u, v)$ is adjacent to $(0, 1)$ 
  thus $(u, v)\not\in \Phi$, which 
  is impossible. We conclude that $\gamma_{t}(\Gamma)>5$. Since vertex $(p-1, p-1)$ is adjacent to vertices
   $(0, 1), (1, p)$ and also vertex $(p-1, 2p-1)$ is adjacent to vertices $(0, 0), (1, p+1)$. Hence
    $T=\{(0,0), (0, 1), (1, p), (1, p+1), (p-1, p-1), (p-1, 2p-1)\}$ is a $\gamma_{t}$-set for $\Gamma$.
    
    Finally $(\alpha, \beta)\neq (1, 1)$. By \cite[Proposition 3.3]{me}, $\gamma(\Gamma)=6$ and 
    $D^{'}=\{(0, 0), (0, 1), (1, 2), (1, 3),\\ (2, 4), (2, 5)\}$ is a $\gamma$-set for $\Gamma$. If $p=3$, then
     vertices $(0, 0), (0, 1), (1, 3)$ are adjacent to vertices $(2, 5), (1, 2), (2, 4)$, respectively and if $p \geq 5 $ 
     then vertices $(0, 0), (1, 3), (2, 4), (0, 1), (1, 2)$ are adjacent to vertices $(1, 3), (2, 4), (0, 1), (1, 2), (2, 5)$, 
     respectively. Thus $D^{'}$ becomes a $\gamma_{t}$-set 
for $\Gamma$. 

Note that both $T$ and $D^{'}$ are two $\gamma_{t}$-sets for $\Gamma$, where $\alpha, \beta\geq 1.$
Therefore $\gamma_{t}(\Gamma)=6$. 

 $2)$ By using a similar argument given in the proof of case $1$, we have $\gamma_{c}(\Gamma)\geq 6$.
 
 Assume first that $p=3$. Then the subgraphs generated by $T$ and $D^{'}$ are disconnected. Since the subgraph generated by $D^{'}$ has exactly 
 three connected components which are induced subgraphs generated by sets $\{(0, 0), (2, 5)\}$, $\{(0, 1), (1, 2)\}$ and $\{(1, 3), (2, 4)\}$, also the subgraph generated by $T$
  has exactly two connected components which are induced subgraphs generated by sets $\{(0, 1), (1, p), (p-1, p-1)\}$ and $\{(0, 0), (1, p+1), (p-1, 2p-1)\}$. 
  We conclude that $\gamma_{c}(\Gamma)\geq 7$. 
 
 Note that vertex $(0, p)$ is adjacent to vertices $(p-1, p-1)$ and $(1, p+1)$. Therefore 
 $C=\{(0, 0), (0, 1), (1, p), (1, p+1), (p-1, p-1), (p-1, 2p-1), (0, p)\}$ is a connected dominating set for 
 $\Gamma$ with minimum cardinality. Therefore $\gamma_{c}(\Gamma)=7$.
 
 Now suppose that $p \geq 5.$ According to the proof of final part of case $1$, we see that $D^{'}$ becomes
  a connected dominating set for $\Gamma$ with minimum cardinality. Therefore in this case  $\gamma_{c}(\Gamma)=6$.
\end{proof}
\end{prop}
\begin{prop}\label{d4}
Let $\Gamma=\Cay(\mathbb{Z}_{p}\times\mathbb{Z}_{p^{\alpha}q^{\beta}}, \Phi)$, where
$p, q\geq 3$ and $\alpha, \beta\geq 1$. Then $diam(\Gamma)=2$.
\begin{proof}
Let $(u, v), (u^{'}, v^{'})\in V(\Gamma)$. Then we have the following three possibilities:

$i)$ $u=u^{'}$ and $v\neq v^{'}$. Hence $d((u, v)$, $(u^{'}, v^{'}))\geq 2$. Let $v$ and $v^{'}$ be
 multiple of $pq$, then $(u-1, pq-1)$ is common neighbor of $(u, v)$ and $(u^{'}, v^{'})$. Let 
 $v, v^{'}\in \varphi_{p^{\alpha}q^{\beta}}$, then $(u-1, pq)$ is adjacent to both $(u, v)$ and $(u^{'}, v^{'})$.
  Let $v$ and $v^{'}$ be multiple of $p$, then $q$ is adjacent to both $v$ and $v^{'}$. Also let $v$ and $v^{'}$ 
  be multiple of $q$, then $p$ is adjacent to both $v$ and $v^{'}$. So $d((u, v)$, $(u^{'}, v^{'}))=2$. Let $v$ is
   multiple of $p$ and $v^{'}$ is multiple of $q$. If  $v$ and $v^{'}$ be both even or odd, then we show that 
   $(u-1, \dfrac{v+v^{'}}{2})$ is a common neighbor of $(u, v)$ and $(u^{'}, v^{'})$. Assume that $v=kp$ and $v^{'}=k^{'}q;~ k, k^{'}\in\mathbb{Z}$. Then $v-\dfrac{v+v^{'}}{2}=\dfrac{v-v^{'}}{2}=\dfrac{kp-k^{'}q}{2}$. 
   Suppose that $\dfrac{kp-k^{'}q}{2}\notin\varphi_{p^{\alpha}q^{\beta}}$ and without loss of generality assume
    $\dfrac{kp-k^{'}q}{2}=k^{''}p;~ k^{''}\in\mathbb{Z}$. Then $kp-k^{'}q=2k^{''}p$ which implies $kp-2k^{''}p=k^{'}q$.
    Hence $(\dfrac{k-2k^{''}}{k^{'}})p=q$, which is impossible, since $q$ is not a multiple of $p$. Hence 
    $\dfrac{kp-k^{'}q}{2}\in \varphi_{p^{\alpha}q^{\beta}}$, and $v$ is adjacent to 
    $\dfrac{v+v^{'}}{2}$. Similarly $v^{'}$ is adjacent to $\dfrac{v+v^{'}}{2}$. If one of either 
     $v$ or $v^{'}$ be odd, then $2(v+v^{'})$ is adjacent to both $v$ and $v^{'}$. 
     Assume that $v=kp$ and even also $v^{'}=k^{'}q$ and odd. 
     Without loss of generality let $2(v+v^{'})-v=v+2v^{'}=k^{''}p$. Then $kp+2k^{'}q=k^{''}p$. 
     This implies $(\dfrac{k^{''}-k}{2k^{'}})p=q$, which is impossible. Thus $v$ is adjacent to $2(v+v^{'})$. 
     Similarly $v^{'}$ is adjacent to $2(v+v^{'})$. Hence
      $d((u, v)$, $(u^{'}, v^{'}))=2$. Let $v$ be
 multiple of $p$ or $q$ or $pq$ and $v^{'}\in \varphi_{p^{\alpha}q^{\beta}}$. 
 Assume that $v$ and $v^{'}$ be both even or odd. If $v-v^{'}\in \varphi_{p^{\alpha}q^{\beta}}$ then it is 
 easy to see that $\dfrac{v+v^{'}}{2}$ is adjacent to both $v$ and $v^{'}$ and if $v-v^{'}\notin \varphi_{p^{\alpha}q^{\beta}}$ then 
 $v-v^{'}$ is adjacent to both $v$ and $v^{'}$. Now suppose that one of either $v$ or $v^{'}$ is odd. If $v$ be 
 multiple of $p$ then $v^{'}q$ is adjacent to both $v$ and $v^{'}$. If $v$ be multiple of $q$ then 
 $v^{'}p$ is adjacent to both $v$ and $v^{'}$. Moreover if $v$ be multiple of $pq$ then $v^{'}(p+q)$ is 
 adjacent to both $v$ and $v^{'}$. Thus $d((u, v)$, $(u^{'}, v^{'}))=2$. Let one of either $v$ or $v^{'}$ is multiple 
 of $p$ or $q$ and other is multiple of $pq$. We know that $+2$ and $-2$ is adjacent to all of the multiple of
  $pq$. Since by proof of \cite[Proposition 3.1]{me}, $\lambda=2$, hence $v$ is adjacent to $+2$ or $-2$ or 
  both of them. So we have a common neighbor between $(u, v)$ and $(u^{'}, v^{'})$. Therefore 
  $d((u, v)$, $(u^{'}, v^{'}))=2$.  
  
  $ii)$  $u\neq u^{'}$ and $v=v^{'}$. In this case the vertex $(u^{''}, v-1)$ where $u^{''}\neq u, u^{'}$,
   is a common neighbor of $(u, v)$ and $(u^{'}, v^{'})$. Thus $d((u, v)$, $(u^{'}, v^{'}))=2$.  
  
  $iii)$ $u\neq u^{'}$ and $v\neq v^{'}$. Hence by $(i)$ and $(ii)$,  $d((u, v)$, $(u^{'}, v^{'}))=2$.
  
  Therefore $diam(\Gamma)=2$.
\end{proof}
\end{prop}
\begin{prop}\label{t3}
Let $\Gamma=\Cay(\mathbb{Z}_{p}\times\mathbb{Z}_{p^{\alpha}q^{\beta}}, \Phi)$, where
$p, q\geq 3$ and $\alpha, \beta\geq 1$. Then $\gamma_{t}(\Gamma)$ and $\gamma_{c}(\Gamma)$ is given by Table~\ref{tab:1}.
\begin{table}
\begin{center}
\caption{$\gamma_{t}(\Cay(\mathbb{Z}_{p}\times \mathbb{Z}_{p^{\alpha}q^{\beta}}, \Phi))=\gamma_{c}(\Cay(\mathbb{Z}_{p}\times \mathbb{Z}_{p^{\alpha}q^{\beta}}, \Phi))$ where $p, q \geq 3$
and $\alpha, \beta\geq 1.$}\label{tab:1}
\begin{tabular}{|c|c|c|}
\hline
$\Gamma$& $\gamma_{t}(\Gamma), \gamma_{c}(\Gamma)$ & Comments \\
\hline
$\Cay(\mathbb{Z}_{p}\times \mathbb{Z}_{p^{\alpha}q^{\beta}}, \Phi)$ & $5$ & one of the prime factors is $3$\\
\hline
$\Cay(\mathbb{Z}_{p}\times \mathbb{Z}_{p^{\alpha}q^{\beta}}, \Phi)$ & $4$ & $p, q \geq 5 $  \\
\hline
\end{tabular}
\end{center}
\end{table}
  \begin{proof}
  Assume first that one of the prime factors is $3$. Let  $(\alpha, \beta)=(1, 1)$. Then by \cite[Proposition 3.1]{me}, 
  $\gamma(\Gamma)=4$ and $D=\{(0, 0), (0, 1), (1, x'), (1, y')\}$ is a $\gamma$-set for
   $\Gamma$, where $x, x'$ and $y, y'$ are consecutive integers in
$\mathbb{Z}_{pq}$, each of which shares a prime factor with $pq$ where $x'$
is a multiple of $p$ and $y'$ is a multiple of $q$. 
  Note that vertices of $D$ are not adjacent to each other. Hence $\gamma_{t}(\Gamma)> 4$. Also 
   $D$ is dominated by $\{(2, 2)\}$. Thus $T=\{(0, 0), (0, 1), (1, x'), (1, y'), (2, 2)\}$ is a
    $\gamma_{t}$-set and $\gamma_{c}$-set
for $\Gamma$.

The next case is where $(\alpha, \beta)\neq (1, 1)$. By \cite[Table 1]{me}, $\gamma(\Gamma)=5$ 
and $D=\{(0, 0), (0, 1), (1, 2)\\, (2, 3), (2, 4)\}$ is a $\gamma$-set for $\Gamma$. 
Vertices $(0, 0), (2, 4), (1, 2), (0, 1)$ are adjacent to vertices $(2, 4), (1,\\ 2), (0, 1), (2, 3)$, 
respectively. Hence $D$ dominates all vertices of $\Gamma$ and the subgraph generated by $D$ is connected. Thus 
$D$ becomes a $\gamma_{t}$-set and $\gamma_{c}$-set
for $\Gamma$. Therefore $\gamma_{t}(\Gamma)=\gamma_{c}(\Gamma)=5$.

Finally assume that $p, q \geq 5$. Then by \cite[Proposition 3.1, Table 1]{me}, $\gamma(\Gamma)=4$ 
and $D=\{(0, 0), (1, 1), (2, 2), (3, 3)\}$ is a $\gamma$-set for $\Gamma$. Since $p, q \geq 5$ then vertices of $D$ 
dominate among themselves. Therefore $\gamma_{t}(\Gamma)=\gamma_{c}(\Gamma)=4$.
  \end{proof}
  \end{prop}
  As an immediate consequence of Lemma~\ref{t1} and Propositions~\ref{t2},~\ref{t3}, 
we have the following theorem.
  \begin{thm}
Let $\Gamma=\Cay(\mathbb{Z}_{p}\times\mathbb{Z}_{p^{\alpha}q^{\beta}}, \Phi)$, where
 $p, q \geq 2$ and $\alpha, \beta\geq 1$. Then
 $\gamma_{t}(\Gamma)$ and $\gamma_{c}(\Gamma)$ is given by Table~\ref{tab:2}.
 \begin{table}
\begin{center}
\caption{$\gamma_{t}(\Cay(\mathbb{Z}_{p}\times \mathbb{Z}_{p^{\alpha}q^{\beta}}, \Phi))$, $\gamma_{c}(\Cay(\mathbb{Z}_{p}\times \mathbb{Z}_{p^{\alpha}q^{\beta}}, \Phi))$ where $\alpha, \beta\geq 1.$}\label{tab:2}
\begin{tabular}{|c|c|c|c|}
\hline
$\Gamma$ & $\gamma_{t}(\Gamma)$ & $\gamma_{c}(\Gamma)$ & Comments \\
\hline
$\Cay(\mathbb{Z}_{2}\times \mathbb{Z}_{2^{\alpha}q^{\beta}}, \Phi)$ & $8$ & does not exist & \\
\hline
$ \Cay(\mathbb{Z}_{p}\times \mathbb{Z}_{2^{\alpha}p^{\beta}}, \Phi)$ & $6$ & $7$ & $p=3$\\
\hline
$ \Cay(\mathbb{Z}_{p}\times \mathbb{Z}_{2^{\alpha}p^{\beta}}, \Phi)$ & $6$ & $6$ & $p\geq 5$\\
\hline
$ \Cay(\mathbb{Z}_{p}\times \mathbb{Z}_{p^{\alpha}q^{\beta}}, \Phi)$ & $5$ & $5$ & one of the prime factors is $3$\\
\hline
$ \Cay(\mathbb{Z}_{p}\times \mathbb{Z}_{p^{\alpha}q^{\beta}}, \Phi)$ & $4$ & $4$ & $p, q\geq 5$\\
\hline
 \end{tabular}
\end{center}
\end{table}
\end{thm}

\begin{ex}
The graph $\Gamma=\Cay(\mathbb{Z}_{2}\times \mathbb{Z}_{2\times 3^{2}},\Phi)$,
which is shown in Figure \ref{fig:4}, is a disconnected graph with two connected components, say $\Gamma_{1}$ and 
$\Gamma_{2}$. Thus  {\it $\gamma_{c}$-set} does not exist for $\Gamma$. 
 In this graph two sets $T_{1}=\{(0, 0), (0, 4), (1, 1), (1, 3)\}$ and
$T_{2}=\{(0, 1), (0, 3), (1, 0), (1, 4)\}$ are {\it $\gamma_{t}$-sets} sets for $\Gamma_{1}$ and $\Gamma_{2}$,
 respectively. Hence $\gamma_{t}(\Gamma)=8$.

\begin{figure}[h]
\begin{center}
\includegraphics[width=7cm , height=6cm]{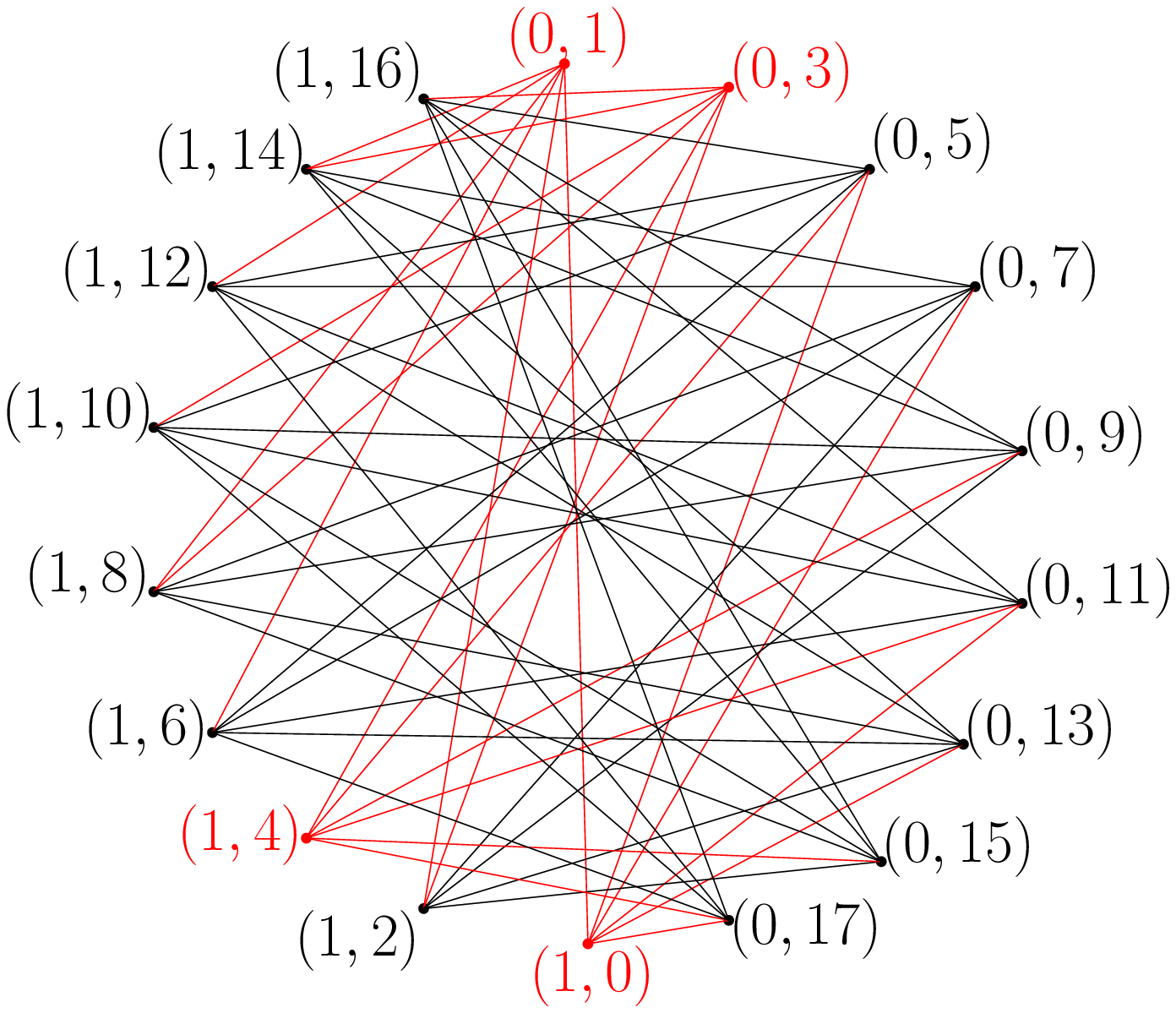} \includegraphics[width=7cm , height=6cm]{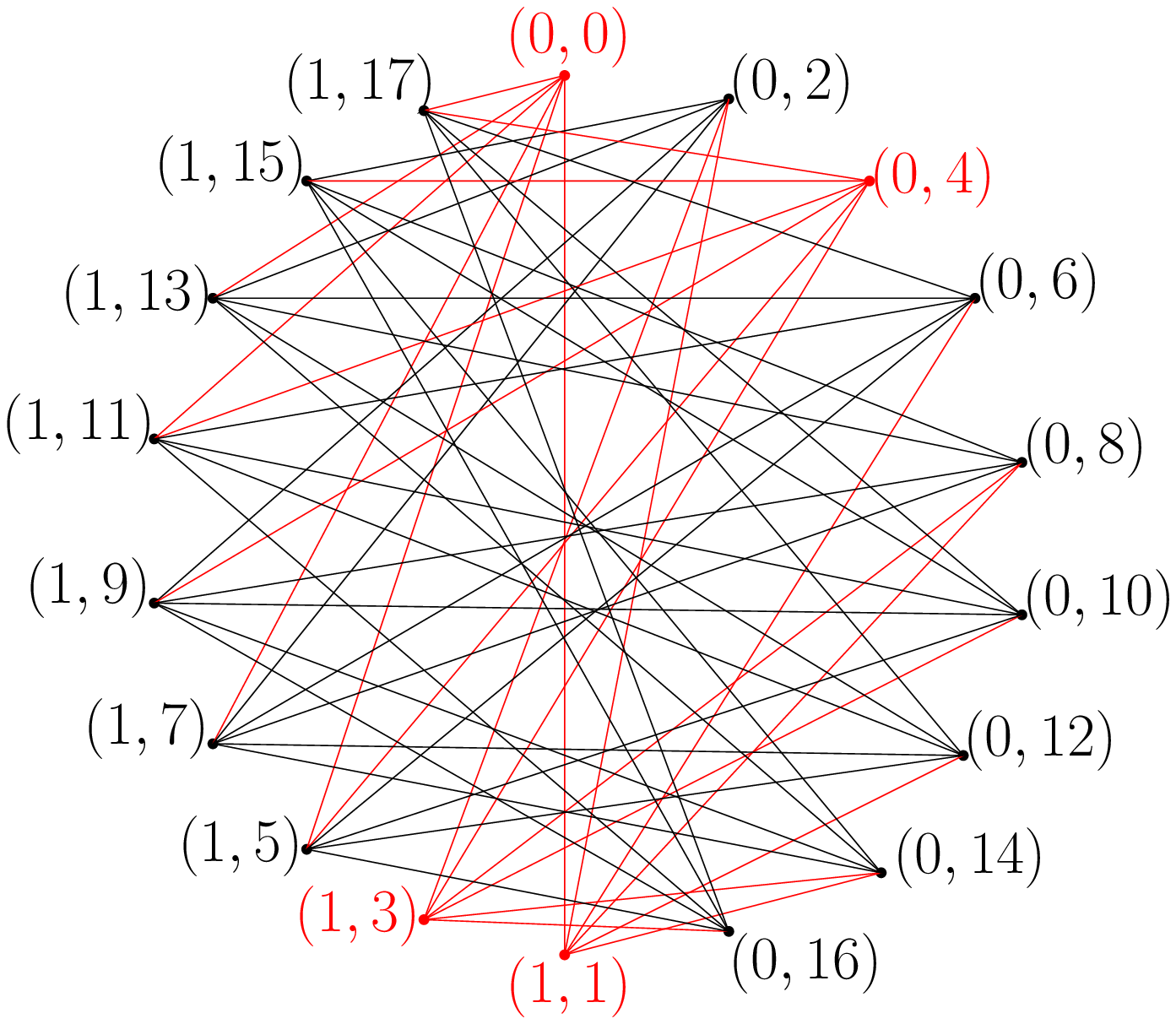}
%\vspace*{-.2cm}
\caption{ Two connected components of $\Gamma=\Cay(\mathbb{Z}_{2}\times \mathbb{Z}_{2\times 3^{2}},\Phi)$, left $\Gamma_{1}$, right $\Gamma_{2}$}\label{fig:4}
\end{center}
\end{figure}
\end{ex}
\begin{ex}
Let $p=3,~q=5$. Then total and connected dominating set of $\Gamma=\Cay(\mathbb{Z}_{3}\times \mathbb{Z}_{15},\Phi)$,
which is shown in Figure \ref{fig:3}, is $\{(0, 0), (0, 1), (1, 6), (1, 10), (2, 2)\}$.

\begin{figure}[h]
\begin{center}
\includegraphics[width=12cm , height=10cm]{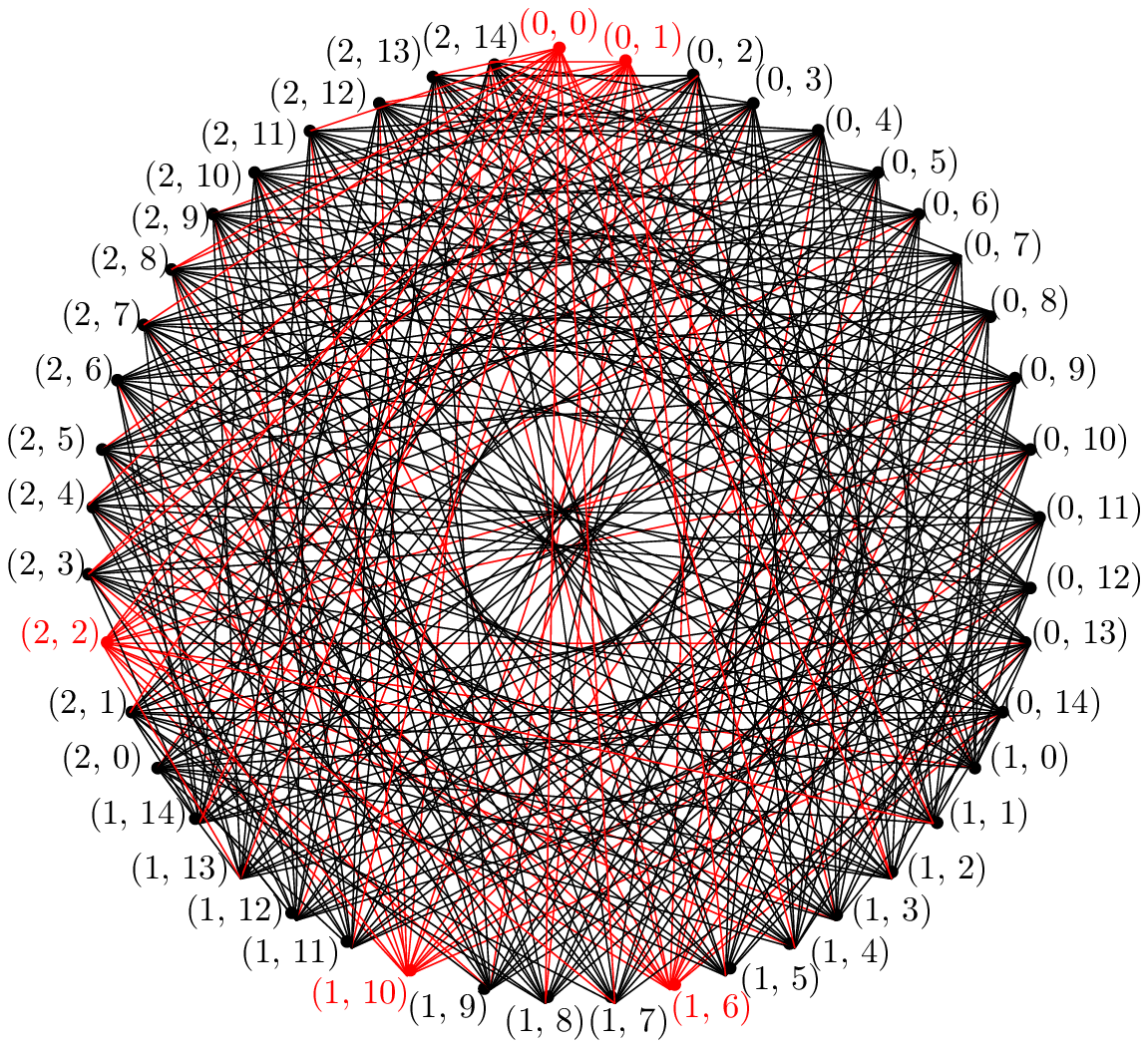}
%\vspace*{-.2cm}
\caption{The graph $\Gamma=\Cay(\mathbb{Z}_{3}\times \mathbb{Z}_{15},\Phi)$ 
and its total dominating set.}\label{fig:3}
\end{center}
\end{figure}
\end{ex}
\section{Total and connected domination number and diameter of 
$\Cay(\mathbb{Z}_{p}\times \mathbb{Z}_{p^{\alpha}q^{\beta}r^{\gamma}},\Phi)$}\label{sec:4}
Let $p, q, r$ be three prime numbers, $\alpha, \beta, \gamma$ positive integers and
$\Phi=\varphi_{p}\times \varphi_{p^{\alpha}q^{\beta}r^{\gamma}}$.
In this section, we obtain the total and connected domination number of
$\Cay(\mathbb{Z}_{p}\times \mathbb{Z}_{p^{\alpha}q^{\beta}r^{\gamma}}, \Phi)$ and we extend the
 results in the previous section for diameter of this graph.
\begin{lem}\label{tdia}
Let $\Gamma=\Cay(\mathbb{Z}_{2}\times \mathbb{Z}_{2^{\alpha}q^{\beta}r^{\gamma}},\Phi)$, where
 $\alpha, \beta, \gamma\geq 1$. Then $diam(\Gamma)=3$.
\begin{proof}
$\Gamma$ is a disconnected graph with two connected
components, say $\Gamma_{1}$ and $\Gamma_{2}$, where
$V(\Gamma_{1})=\{(1, v)| v~is~odd\}\cup \{(0, v)| v~is~even\}$ and
$V(\Gamma_{2})=\{(0, v)| v~is~odd\}\cup \{(1, v)| v~is~even\}$.

Let $(u, v), (u^{'}, v^{'})\in V(\Gamma_{1})$. Then we have the following two possibilities:

$i)$ $u=u^{'}, v\neq v^{'}$. Since $u=u^{'}$ hence $d((u, v), (u^{'}, v^{'}))\geq 2$. Now by Table~\ref{tab:dia1} 
we show that $d((u, v), (u^{'}, v^{'}))=2$. In this table, when $v, v^{'}$ are odd we have $u=u^{'}=1,~ u^{''}=0$ and when  $v, v^{'}$ are even we have $u=u^{'}=0,~u^{''}=1$. We prove the rows $6, 8$ of the table and the rest is similarly proven.  

Let $v, v^{'}$ are odd and $v=kq,~ v^{'}=k^{'}qr,~k, k^{'}\in\mathbb{Z}$. If $\dfrac{v+v^{'}}{q}$ be non-multiple of $q$ then we show that $\dfrac{v+v^{'}}{q}$ is adjacent to both $v$ and $v^{'}$.

Let $k^{''}\in\mathbb{Z}$. If $v-\dfrac{v+v^{'}}{q}=2k^{''}$, then $k(q-1)-k^{'}r=2k^{''}$. This implies $k=\dfrac{2k^{''}+k^{'}r}{q-1}$. Since $2k^{''}+k^{'}r$ is odd and $q-1$ is even hence $k$ is non-integer, which is impossible. If $v-\dfrac{v+v^{'}}{q}=k^{''}q$, then $\dfrac{v+v^{'}}{q}=(k-k^{''})q$, which is inaccurate because $\dfrac{v+v^{'}}{q}$ is non-multiple of $q$. Moreover if $v-\dfrac{v+v^{'}}{q}=k^{''}r$, then $k=(\dfrac{k^{'}+k^{''}}{q-1})r$. But we know that $k$ is non-multiple of $r$. So $v-\dfrac{v+v^{'}}{q}\in\varphi_{2^{\alpha}q^{\beta}r^{\gamma}}$ and similarly $v^{'}$ is adjacent to $\dfrac{v+v^{'}}{q}$. Since $u^{''}$ is adjacent to $u, u^{'}$ thus $(u^{''}, \dfrac{v+v^{'}}{q})$ is common neighbor between $(u, v), (u^{'}, v^{'})$. Similarly it is easy to see that if $\dfrac{v+v^{'}}{q}$ be multiple of $q$ then $(u^{''}, \dfrac{v+v^{'}}{q}+2r)$ is adjacent to both $(u, v)$ and $(u^{'}, v^{'})$.

Let $v \in\varphi_{2^{\alpha}q^{\beta}r^{\gamma}}, v^{'}=kq$ is odd and $~k, k^{''}\in\mathbb{Z}$. If $v-(v+v^{'})r=2k^{''}$, then $v=2k^{''}+(v+v^{'})r$. Hence $v$ is even, which is inaccurate. Also if $v-(v+v^{'})r=k^{''}q$, then $v^{'}=(\dfrac{k(1-r)-k^{''}}{r})q$ and if $v-(v+v^{'})r=k^{''}r$, $v=(k^{''}+v+v^{'})r$, which are impossible. Hence $v$ is adjacent to $(v+v^{'})r$. Similarly it is easy to see that $v^{'}$ is adjacent to $(v+v^{'})r$. Therefore $(u^{''}, (v+v^{'})r)$ is adjacent to both $(u, v)$ and $(u^{'}, v^{'})$. 

$ii)$ $u\neq u^{'}, v\neq v^{'}$. If $v$ be adjacent to $v^{'}$, then $d((u, v), (u^{'}, v^{'}))=1$. 
  Suppose that $v$ be non-adjacent to $v^{'}$, since $u\neq u^{'}$ and $u, u^{'}\in \mathbb{Z}_{2}$, hence
   we have no common neighbor between $(u, v)$ and  $(u^{'}, v^{'})$. 
   This implies that $d((u, v), (u^{'}, v^{'}))\geq 3$. 
   Without loss of generality assume that $u=0$ and $u^{'}=1$. Now by Table~\ref{tab:dia2} 
   we show that $d((u, v), (u^{'}, v^{'}))=3$. In this table $u, u^{'''}=0$ and also $u^{'}, u^{''}=1$.
    Now we prove the fifth row and the rest is similarly proven. Let $v=2kr;~k\in\mathbb{Z}$ and $v^{'}\in\varphi_{2^{\alpha}q^{\beta}r^{\gamma}}$. Clearly $u, u^{'''}$ are adjacent to $u^{'}, u^{''}$. 

First we show that $v$ is adjacent to $q$. Let $k^{''}\in\mathbb{Z}$.
$$ \text{If}~ v-q=2k^{''}, ~\text{then}~ q=2(kr-k^{''}) $$
$$ \text{If}~ v-q=k^{''}q, ~\text{then}~ r=(\dfrac{k^{''}+1}{2k})q$$
$$ \text{If}~ v-q=k^{''}r, ~\text{then}~ q=(2k-k^{''})r.$$
In all three cases, we came across a contradiction. So $v-q\in\varphi_{2^{\alpha}q^{\beta}r^{\gamma}}$.

Next we prove that $q$ is adjacent to $(q+v^{'})r$.
$$ \text{If}~ (q+v^{'})r-q=2k^{''}, ~\text{then}~ k^{''}=\dfrac{(q+v^{'})r-q}{2}$$
$$  \text{If}~ (q+v^{'})r-q=k^{''}q, ~\text{then}~ v^{'}=(\dfrac{k^{''}-r+1}{r})q$$
$$  \text{If}~ (q+v^{'})r-q=k^{''}r, ~\text{then}~ q=(q+v^{'}-k^{''})r.$$
Which are impossible. Since $k^{''}$ is integer and $v^{'}\in\varphi_{2^{\alpha}q^{\beta}r^{\gamma}}$ 
and also $q$ is non-integer of $r$. 

Finally we show that $(q+v^{'})r$ is adjacent to $v^{'}$.
$$ \text{If}~(q+v^{'})r-v^{'}=2k^{''}, ~\text{then}~ k^{''}=\dfrac{(q+v^{'})r-v^{'}}{2}$$
$$ \text{If}~(q+v^{'})r-v^{'}=k^{''}q, ~\text{then}~ v^{'}=(\dfrac{k^{''}-r}{r-1})q$$
$$ \text{If}~(q+v^{'})r-v^{'}=k^{''}r, ~\text{then}~ v^{'}=(q+v^{'}-k^{''})r.$$
Again which are impossible. This implies 
that $(u, v)(u^{''}, q)(u^{'''}, (q+v^{'})r)(u^{'}, v^{'})$ is shortest path between $(u, v)$ 
and $(u^{'}, v^{'})$. Thus $diam(\Gamma_{1})=3$ and similarly $diam(\Gamma_{2})=3$. 
Therefore $diam(\Gamma)=3$.
 \begin{table}
\begin{center}
\caption{common neighbor between $(u, v), (u^{'}, v^{'})$ in $\Gamma=\Cay(\mathbb{Z}_{2}\times \mathbb{Z}_{2^{\alpha}q^{\beta}r^{\gamma}},\Phi)$}\label{tab:dia1}
\begin{tabular}{|c|c|c|}
\hline
$u=u^{'}, v\neq v^{'}$ & common neighbor & Comments \\
\hline
 $v, v^{'}\in \varphi_{2^{\alpha}q^{\beta}r^{\gamma}}$& $(u^{''}, 2qr)$ &\\
\hline
$v, v^{'}$ are odd and multiple of $q$ & $(u^{''}, 2r)$ &\\
\hline
$v, v^{'}$ are odd and multiple of $r$ & $(u^{''}, 2q)$ & \\
\hline
$v, v^{'}$ are odd and multiple of $qr$ & $(u^{''}, 2)$ &\\

\hline
$v$ is odd and multiple of $q$ & $(u^{''}, \dfrac{v+v^{'}}{2})$& if $\dfrac{v+v^{'}}{2}$ be even\\
and $v^{'}$ is odd and multiple of $r$ & ~$(u^{''}, \dfrac{v+v^{'}}{2}+qr)$ & if $\dfrac{v+v^{'}}{2}$ be odd\\
\hline
$v$ is odd and multiple of $q$ & $(u^{''}, \dfrac{v+v^{'}}{q})$& if $\dfrac{v+v^{'}}{q}\neq kq,~k\in\mathbb{Z}$\\
and $v^{'}$ is odd and multiple of $qr$ & ~$(u^{''}, \dfrac{v+v^{'}}{q}+2r)$ & if $\dfrac{v+v^{'}}{q}=k^{'}q,~k^{'}\in\mathbb{Z}$\\
\hline
$v$ is odd and multiple of $r$ & $(u^{''}, \dfrac{v+v^{'}}{r})$& if $\dfrac{v+v^{'}}{r}\neq kr,~k\in \mathbb{Z}$\\
and $v^{'}$ is odd and multiple of $qr$ & ~$(u^{''}, \dfrac{v+v^{'}}{r}+2q)$ & if $\dfrac{v+v^{'}}{r}=k^{'}r,~k^{'}\in\mathbb{Z}$\\
\hline
$v\in\varphi_{2^{\alpha}q^{\beta}r^{\gamma}}$ and $v^{'}$ is odd and multiple of $q$ & $(u^{''}, (v+v^{'})r)$&\\
\hline
$v\in\varphi_{2^{\alpha}q^{\beta}r^{\gamma}}$ and $v^{'}$ is odd and multiple of $r$ & $(u^{''}, (v+v^{'})q)$&\\
\hline
$v\in\varphi_{2^{\alpha}q^{\beta}r^{\gamma}}$ and $v^{'}$ is odd multiple of $qr$ & $(u^{''}, (v+v^{'})2)$&\\
\hline
$v, v^{'}$ are even and multiple of $r$ & $(u^{''}, q)$ &\\
\hline
$v, v^{'}$ are even and multiple of $q$ & $(u^{''}, r)$ & \\
\hline
$v, v^{'}$ are even and non-multiple of $q$ and $r$ & $(u^{''}, qr)$ & \\
\hline
$v, v^{'}$ are even and multiple of $2qr$ & $(u^{''}, 2qr-1)$ &\\
\hline
$v$ is even and multiple of $q$ & $(u^{''}, \dfrac{v+v^{'}}{2}+qr)$& if $\dfrac{v+v^{'}}{2}$ be even\\
and $v^{'}$ is even and multiple of $r$ & $(u^{''}, \dfrac{v+v^{'}}{2})$~~~~~ & if $\dfrac{v+v^{'}}{2}$ be odd\\
\hline
$v$ is even and multiple of $qr$ & $(u^{''}, \dfrac{v^{'}}{2})$& if $\dfrac{v^{'}}{2}\in \varphi_{2^{\alpha}q^{\beta}r^{\gamma}}$\\
and $v^{'}$ is even and non-multiple of & $(u^{''}, \dfrac{v^{'}}{2}+qr)$ & if $\dfrac{v^{'}}{2}\notin \varphi_{2^{\alpha}q^{\beta}r^{\gamma}}$\\
$q$ and $r$ &&\\
\hline
$v$ is even and multiple of $qr$ & $(u^{''}, \dfrac{v^{'}}{q}+qr)$& if $\dfrac{v^{'}}{q}\neq kq, ~k\in\mathbb{Z}$\\
and $v^{'}$ is even and multiple of $q$ & $(u^{''}, \dfrac{v^{'}}{q}+r)$ & if $\dfrac{v^{'}}{q}=k^{'}q,~k^{'}\in\mathbb{Z}$\\
\hline
$v$ is even and multiple of $qr$ & $(u^{''}, \dfrac{v^{'}}{r}+qr)$& if $\dfrac{v^{'}}{r}\neq kr,~k\in\mathbb{Z}$\\
and $v^{'}$ is even and multiple of $r$ & $(u^{''}, \dfrac{v^{'}}{r}+q)$ & if $\dfrac{v^{'}}{r}=k^{'}r,~k^{'}\in\mathbb{Z}$\\
\hline
\end{tabular}
\end{center}
\end{table}
\begin{table}
\begin{center}
\caption{shortest path between $(u, v), (u^{'}, v^{'})$ in $\Gamma=\Cay(\mathbb{Z}_{2}\times \mathbb{Z}_{2^{\alpha}q^{\beta}r^{\gamma}},\Phi)$}\label{tab:dia2}
\begin{tabular}{|c|c|c|}
\hline
$u\neq u^{'}, v\neq v^{'}$  & shortest path between $(u, v), (u^{'}, v^{'})$ & Comments \\
\hline
 $v=2kqr,k\in\mathbb{Z},v^{'}\in \varphi_{2^{\alpha}q^{\beta}r^{\gamma}}$& $(u, v)(u^{'}, v^{'})$  &$d((u, v), (u^{'}, v^{'}))=1$\\
\hline
$v$ is multiple of $2qr$ and & $(u, v)(u^{''}, 1)(u^{'''}, (1+v^{'})r)(u^{'}, v^{'})$ &\\
$v^{'}$ is odd and multiple of $q$&&\\
\hline
$v$ is multiple of $2qr$ and & $(u, v)(u^{''}, 1)(u^{'''}, (1+v^{'})q)(u^{'}, v^{'})$ &\\
$v^{'}$ is odd and multiple of $r$&&\\
\hline
$v$ is multiple of $2qr$ and & $(u, v)(u^{''}, 1)(u^{'''}, (1+v^{'})2)(u^{'}, v^{'})$ &\\
$v^{'}$ is odd and multiple of $qr$&&\\
\hline
 $v=2kr,k\in\mathbb{Z},v^{'}\in \varphi_{2^{\alpha}q^{\beta}r^{\gamma}}$& $(u, v)(u^{''}, q)(u^{'''}, (q+v^{'})r)(u^{'}, v^{'})$  &\\
\hline
$v$ is multiple of $2r$ and & ~~~$(u, v)(u^{''}, q)(u^{'''}, \dfrac{q+v^{'}}{2})(u^{'}, v^{'})$ & if $\dfrac{q+v^{'}}{2}$ be even\\
$v^{'}$ is odd and multiple of $r$&~$(u, v)(u^{''}, q)(u^{'''}, \dfrac{q+v^{'}}{2}+qr)(u^{'}, v^{'})$&if $\dfrac{q+v^{'}}{2}$ be odd\\
\hline
$v$ is multiple of $2r$ and & $(u, v)(u^{''}, q)(u^{'''}, 2r)(u^{'}, v^{'})$ &\\
$v^{'}$ is odd and multiple of $q$&&\\
\hline
$v$ is multiple of $2r$ and & ~~~$(u, v)(u^{''}, q)(u^{'''}, \dfrac{q+v^{'}}{q})(u^{'}, v^{'})$ & if $\dfrac{q+v^{'}}{q}\neq kq,k\in\mathbb{Z}$\\
$v^{'}$ is odd and multiple of $qr$&$(u, v)(u^{''}, q)(u^{'''}, \dfrac{q+v^{'}}{q}+2r)(u^{'}, v^{'})$&if $\dfrac{q+v^{'}}{q}=k^{'}q,k^{'}\in\mathbb{Z}$\\
\hline
 $v=2kq,k\in\mathbb{Z},v^{'}\in \varphi_{2^{\alpha}q^{\beta}r^{\gamma}}$& $(u, v)(u^{''}, r)(u^{'''}, (r+v^{'})q)(u^{'}, v^{'})$  &\\
\hline
$v$ is multiple of $2q$ and & ~~~$(u, v)(u^{''}, r)(u^{'''}, \dfrac{r+v^{'}}{2})(u^{'}, v^{'})$ & if $\dfrac{r+v^{'}}{2}$ be even\\
$v^{'}$ is odd and multiple of $q$&~$(u, v)(u^{''}, r)(u^{'''}, \dfrac{r+v^{'}}{2}+qr)(u^{'}, v^{'})$&if $\dfrac{r+v^{'}}{2}$ be odd\\
\hline
$v$ is multiple of $2q$ and & $(u, v)(u^{''}, r)(u^{'''}, 2q)(u^{'}, v^{'})$ &\\
$v^{'}$ is odd and multiple of $r$&&\\
\hline
$v$ is multiple of $2q$ and & ~~~$(u, v)(u^{''}, r)(u^{'''}, \dfrac{r+v^{'}}{r})(u^{'}, v^{'})$ & if $\dfrac{r+v^{'}}{r}\neq kr,k\in\mathbb{Z}$\\
$v^{'}$ is odd and multiple of $qr$&~$(u, v)(u^{''}, r)(u^{'''}, \dfrac{r+v^{'}}{r}+2q)(u^{'}, v^{'})$&if $\dfrac{r+v^{'}}{r}=k^{'}r,k^{'}\in\mathbb{Z}$\\
\hline
 $v=2k,k\in\mathbb{Z},v^{'}\in \varphi_{2^{\alpha}q^{\beta}r^{\gamma}}$& $(u, v)(u^{''}, qr)(u^{'''}, (qr+v^{'})2)(u^{'}, v^{'})$  &\\
\hline
$v$ is multiple of $2$ and & ~~~$(u, v)(u^{''}, qr)(u^{'''}, \dfrac{qr+v^{'}}{q})(u^{'}, v^{'})$ & if $\dfrac{qr+v^{'}}{q}\neq kq,k\in\mathbb{Z}$\\
$v^{'}$ is odd and multiple of $q$&$(u, v)(u^{''}, qr)(u^{'''}, \dfrac{qr+v^{'}}{q}+2r)(u^{'}, v^{'})$&if $\dfrac{qr+v^{'}}{q}=k^{'}q,k^{'}\in\mathbb{Z}$\\
\hline
$v$ is multiple of $2$ and & ~~~$(u, v)(u^{''}, qr)(u^{'''}, \dfrac{qr+v^{'}}{r})(u^{'}, v^{'})$ & if $\dfrac{qr+v^{'}}{r}\neq kr,k\in\mathbb{Z}$\\
$v^{'}$ is odd and multiple of $r$&$(u, v)(u^{''}, qr)(u^{'''}, \dfrac{qr+v^{'}}{r}+2q)(u^{'}, v^{'})$ &if $\dfrac{qr+v^{'}}{r}=k^{'}r,k^{'}\in\mathbb{Z}$\\
\hline
$v$ is multiple of $2$ and & $(u, v)(u^{''}, qr)(u^{'''}, 2)(u^{'}, v^{'})$ &\\
$v^{'}$ is odd and multiple of $qr$&&\\
\hline
\end{tabular}
\end{center}
\end{table}
\end{proof}
\end{lem}
\begin{lem}\label{t4}
Let $\Gamma=\Cay(\mathbb{Z}_{2}\times \mathbb{Z}_{2^{\alpha}q^{\beta}r^{\gamma}},\Phi)$, where
 $\alpha, \beta, \gamma\geq 1$. Then $\gamma_{c}(\Gamma)$ does not exist and $\gamma_{t}(\Gamma)=12$.
\begin{proof}
Clearly $\Gamma$ is a disconnected graph with two connected components say $\Gamma_{1}$ and $\Gamma_{2}$.
Let $V_{1}=V(\Gamma_{1})$ and $V_{2}=V(\Gamma_{2})$. Then 
$V_{1}=\{(1, v)| v\text{~is~odd}\}\cup \{(0, v)| v\text{~is~even}\}$
and $V_{2}=\{(0, v)| v\text{~is~odd}\}\cup \{(1, v)| v\text{~is~even}\}.$ Hence by the definition of 
connected dominating set, {\it $\gamma_{c}$-set} does not exist for $\Gamma$.

Let $(\alpha, \beta, \gamma)=(1, 1, 1)$. By \cite[Lemma 4.1]{me}, $\gamma(\Gamma)=8$ and 
$D_{1}=\{(0, 0), (0, 2), (1, x_{4}), (1, x_{4}')\},\\~D_{2}=\{(0, 1), (0, 3), (1, x_{5}), (1, x_{5}')\}$
are minimal dominating sets for $\Gamma_{1}$ and $\Gamma_{2}$ respectively, where 
$X_{i}^{5}=\{x_{1}, x_{2}, x_{3}, x_{4}, x_{5}\}$ and $X_{j}^{5}=\{x_{1}', x_{2}', x_{3}', x_{4}', x_{5}'\}$ 
are consecutive integers in
$\mathbb{Z}_{2qr}$, each of which shares a prime factor with $2qr$. Since vertices of $D_{1}$ are
 not adjacent to each other, we conclude that $\gamma_{t}(\Gamma_{1})>4$. On the other hand
  it is clear that $D_{1}$ is not dominated by one vertex. Hence $\gamma_{t}(\Gamma_{1})>5$. 
  Vertex $(1, 1)$ is adjacent to vertices $(0, 0), (0, 2)$ and vertex $(0, 4)$ is adjacent to vertices 
  $(1, x_{4}), (1, x_{4}')$. Thus $T_{1}=\{(0, 0), (0, 2), (0, 4), (1, 1), (1, x_{4}), (1, x_{4}'\\)\}$ dominates all 
  vertices of $\Gamma_{1}$. Similarly $T_{2}=\{(0, 1), (0, 3), (0, 5), (1, 2), (1, x_{5}), (1, x_{5}')\}$
   dominates all vertices of $\Gamma_{2}$. Hence $\gamma_{t}(\Gamma)=12$.

Let $(\alpha, \beta, \gamma)\neq (1, 1, 1)$. By \cite[Lemma 4.3]{me}, $\gamma(\Gamma)=12$ and 
$D_{1}=\{(0, 0), (0, 2), (0, 4), (1, 1),\\ (1, 3), (1, 5)\},~ D_{2}=\{(0, 1), (0, 3), (0, 5), (1, 0), (1, 2), (1, 4)\}$ 
are minimal dominating sets for $\Gamma_{1},~\Gamma_{2}$, respectively. Vertex $(1, 1)$ is adjacent to vertices 
$(0, 0), (0, 2)$ and vertex $(0, 4)$ is adjacent to vertices 
  $(1, 3), (1, 5)$. Thus $D_{1}$ becomes a $\gamma_{t}$-set 
for $\Gamma_{1}$. Similarly $D_{2}$ becomes a $\gamma_{t}$-set 
for $\Gamma_{2}$. Therefore $\gamma_{t}(\Gamma)=12$.

\end{proof}
\end{lem}
\begin{prop}\label{d6}
Let $\Gamma=\Cay(\mathbb{Z}_{p}\times \mathbb{Z}_{2^{\alpha}p^{\beta}r^{\gamma}},\Phi)$,
 where $\alpha, \beta, \gamma\geq 1$. Then $diam(\Gamma)=3$.
\begin{proof}
We proceed along the lines of Theorem~\ref{tdia}, and $q:=p$. Let $(u, v), (u^{'}, v^{'})$ are arbitrary vertices of $\Gamma$. 
Then we have following three possibilities:

$i)$ $u=u^{'}, v\neq v^{'}$. We know that $d((u, v), (u^{'}, v^{'}))\geq 2$. Assume that $v$ and $v^{'}$ are 
both even or odd. Thus by case $(i)$ of Theorem~\ref{tdia}, we have $d((u, v), (u^{'}, v^{'}))=2$. Suppose that 
one of either $v$ or $v^{'}$ is odd. Hence we have no path of length $2$ between $(u, v), (u^{'}, v^{'})$. 
Now we show that $d((u, v), (u^{'}, v^{'}))=3$. Without loss of generality assume that $v$ is even and $v^{'}$ is odd. If $v$ 
be multiple of $2pr$ and $v^{'}$ be multiple of $pr$, then $(u, v)(u^{''}, pr-2)(u^{'''}, pr-1)(u^{'}, v^{'})$ is a path of length 
$3$ between $(u, v)$ and $(u^{'}, v^{'})$, where $u=u^{'}\neq u^{''}\neq u^{'''}$. If $v$ 
be multiple of $2pr$ and $v^{'}\in \varphi_{2^{\alpha}p^{\beta}r^{\gamma}}$, note that $v$ and $v^{'}$ are adjacent, then 
$(u, v)(u^{''}, v^{'})(u^{'''}, v)(u^{'}, v^{'})$ is a shortest path. 
For other cases of $v$ and $v^{'}$ we are using of Table~\ref{tab:dia2}, where $u=u^{'}\neq u^{''}\neq u^{'''}$. 

$ii)$ $u\neq u^{'}, v=v^{'}$. In this case vertex $(u^{''}, v-1)$ is a common neighbor between $(u, v)$ and $(u^{'}, v^{'})$, 
where $u^{''}\neq u, u^{'}$. Thus $d((u, v), (u^{'}, v^{'}))=2$.

$iii)$ $u\neq u^{'}, v\neq v^{'}$. Let $v$ and $v^{'}$ be both even or odd. Then by Table~\ref{tab:dia1}, 
where $u^{''}\neq u, u^{'}$, we see that $d((u, v), (u^{'}, v^{'}))=2$. Let one of either $v$ or $v^{'}$ be even 
and other be odd. Then by Table~\ref{tab:dia2}, where $u=u^{'''}$ and $u^{'}=u^{''}$, we see that 
$d((u, v), (u^{'}, v^{'}))=3$. Therefore $diam(\Gamma)=3$.  
\end{proof}
\end{prop}
\begin{prop}\label{t5}
Let $\Gamma=\Cay(\mathbb{Z}_{p}\times \mathbb{Z}_{2^{\alpha}p^{\beta}r^{\gamma}},\Phi)$,
 where $\alpha, \beta, \gamma\geq 1$. Then $\gamma_{t}(\Gamma)$ and $\gamma_{c}(\Gamma)$ is given by Table~\ref{tab:3}.
\begin{table}
\begin{center}
\caption{$\gamma_{t}(\Cay(\mathbb{Z}_{p}\times \mathbb{Z}_{2^{\alpha}p^{\beta}r^{\gamma}}, \Phi))$ and $\gamma_{c}(\Cay(\mathbb{Z}_{p}\times \mathbb{Z}_{2^{\alpha}p^{\beta}r^{\gamma}}, \Phi))$}
\label{tab:3}
\begin{tabular}{|c|c|c|c|}

\hline
$\Gamma$& $\gamma_{t}(\Gamma)$ & $\gamma_{c}(\Gamma)$ & Comments \\
\hline
$\Cay(\mathbb{Z}_{p}\times \mathbb{Z}_{2^{\alpha}p^{\beta}r^{\gamma}}, \Phi)$ & $10$ & $12$ & one of the prime factors is $3$\\
&&& $p=3$\\
\hline
$\Cay(\mathbb{Z}_{p}\times \mathbb{Z}_{2^{\alpha}p^{\beta}r^{\gamma}}, \Phi)$ & $10$ & $10$ & one of the prime factors is $3$\\
&&& $p\geq 5 $  \\
\hline
$\Cay(\mathbb{Z}_{p}\times \mathbb{Z}_{2^{\alpha}p^{\beta}r^{\gamma}}, \Phi)$ & $8$ & $8$ & $p, r \geq 5 $  \\
\hline
\end{tabular}
\end{center}
\end{table}
\begin{proof}
Assume first that one of the prime factors is $3$. In this case if $(\alpha, \beta, \gamma)=(1, 1, 1)$ 
then by \cite[Lemma 4.2]{me}, $\gamma(\Gamma)=8$ and 
$D=\{(0, 0), (0, 1), (0, 2), (0, 3), (1, x_{4}), (1, x_{4}^{'}), (1, x_{5}), (1, x_{5}^{'})\}$ is a $\gamma$-set
 for $\Gamma$. Vertices of $D$ are not adjacent to each other. Hence $\gamma_{t}(\Gamma)>8$.
  Note that $D$ is not dominated by one vertex, since every vertex $(u, v)\in V$, where $v$ is an odd (even) integer,
   is not adjacent to the vertex $(u^{'}, v^{'})$, where $v^{'}$ is an odd (even) integer. This implies that 
   $\gamma_{t}(\Gamma)>9$. Now we take another dominating set with cardinality $10$. By \cite[Proposition 4.4]{me},
    we have $D^{'}=\{(0, 0), (0, 1),(0, 2), (0, 3), (1, 4), (1, 5), (2, 6), (2, 7), (2, 8), (2, 9)\}$ is a dominating set of 
    $\Gamma$. Let other prime factor be $5$ then vertices $(0, 0), (0, 1),(0, 2), (0, 3), (1, 5)$ are adjacent to vertices 
    $(2, 7), (2, 8), (2, 9), (1, 4), (2, 6)$, respectively. Also  let other prime factor be $\geq 7$ then vertices 
    $(0, 0), (0, 1),(0, 2), (0, 3), (1, 4)$ are adjacent to vertices $ (1, 5), (2, 6), (2, 7), (2, 8), (2, 9)$, respectively. 
    Hence $D^{'}$ becomes a $\gamma_{t}$-set 
for $\Gamma$. Therefore $\gamma_{t}(\Gamma)=10$.

Let $(\alpha, \beta, \gamma)\neq (1, 1, 1)$. By  \cite[Proposition 4.4]{me}, $\gamma(\Gamma)=10$. 
By previous paragraph, $\gamma_{t}(\Gamma)=10$.

Now we find the connected domination number of $\Gamma$ where one of the prime factors is $3$. 
By above discussion $\gamma_{c}(\Gamma)>9$. We use again from $D^{'}$. 

Let $p=3$. Without loss of generality assume that $r=5$. Then the subgraph generated by $D^{'}$ has exactly five connected components which are induced the subgraphs generated by sets
$\{(0, 0), (2, 7)\}$, $\{(0, 1), (2, 8)\}$, $\{(0, 2), (2, 9)\}$, $\{(0, 3), (1, 4)\}$ and $\{(1, 5), (2, 6)\}$. Hence
 $\gamma_{c}(\Gamma)>10$. Let a vertex say $(u, v)\in V(\Gamma)$, where $v$ is an odd integer, dominates all 
vertices $(0, 0), (2, 8), (0, 2), (1, 4), (2, 6) $. Since $u\in \mathbb{Z}_{3}$, it is impossible. 
This implies that $\gamma_{c}(\Gamma)>11$. Next consider another dominating with cardinality $12$.

Let $A=\{(1, 1), (2, 2), (1, 4), (2, 5), (1, 7), (2, 8), (1, 10), (2, 11)\}$, 
$B=\{(0, 0), (2, 2), (0, 3), (2, 5)\\, (0, 6), (2, 8), (0, 9), (2, 11)\}$ and 
$C=\{(0, 0), (1, 1), (0, 3), (1, 4), (0, 6), (1, 7), (0, 9), (1, 10)\}$. Then $A$, $B$ and $C$ dominate 
$\{(0, v)| v \in  \mathbb{Z}_{2^{\alpha}3^{\beta}r^{\gamma}}\}$, 
$\{(1, v)| v \in  \mathbb{Z}_{2^{\alpha}3^{\beta}r^{\gamma}}\}$ and
 $\{(2, v)| v \in  \mathbb{Z}_{2^{\alpha}3^{\beta}r^{\gamma}}\}$ respectively.
 
  Thus  
 $D^{''}=\{(0, 0), (1, 1), (2, 2), (0, 3), (1, 4), (2, 5), (0, 6), (1, 7), (2, 8), (0, 9), (1, 10), (2, 11) \}$
  is a dominating set for $\Gamma$. Both vertices next to each other in $D^{''}$ are adjacent.
   Hence the subgraph generated by $D^{''}$ is connected. Therefore $\gamma_{c}(\Gamma)=12$.  
   
Let $p \geq 5 $. Then $D^{'''}=\{(0, 0), (1, 1), (2, 2), (3, 3), (4, 4), (0, 5), (1, 6), (2, 7), (3, 8), (4, 9)\}$ 
is a dominating set for $\Gamma$. Both vertices next to each other in $D^{'''}$ are adjacent. 
Thus the subgraph generated by $D^{'''}$ is connected. Therefore $\gamma_{c}(\Gamma)=10$.    
   
Finally assume that $p, r \geq 5 $. By \cite[Lemma 4.2, Proposition 4.4]{me}, $\gamma(\Gamma)=8$ 
and by using a proof of proposition 4.4, we know that $D^{''''}=\{(0, 0), (1, 1), (2, 2), (3, 3), (4, 4), (2, 5), (1, 6), (0, 7)\}$
 is a $\gamma$-set for $\Gamma$, where $\alpha, \beta, \gamma\geq 1$. Both vertices next to each other in $D^{''''}$ are adjacent. 
Hence $D^{''''}$ is a $\gamma_{t}$-set and $\gamma_{c}$-set 
for $\Gamma$. Therefore $\gamma_{t}(\Gamma)=\gamma_{c}(\Gamma)=8$.
\end{proof}
\end{prop}
\begin{prop}\label{d7}
Let $\Gamma=\Cay(\mathbb{Z}_{p}\times\mathbb{Z}_{p^{\alpha}q^{\beta}r^{\gamma}}, \Phi)$,
where $p, q, r \geq 3 $ and $\alpha, \beta, \gamma\geq 1$.
Then $diam(\Gamma)=2$.
\begin{proof}
Let $(u, v), (u^{'}, v^{'})$ are arbitrary vertices of $\Gamma$. Then we have following three possibilities:

$i)$ $u=u^{'}, v\neq v^{'}$. By Table~\ref{tab:dia3}, we show that $d((u, v), (u^{'}, v^{'}))=2$. In this table $u^{''}\neq u$.

  $ii)$  $u\neq u^{'}$ and $v=v^{'}$. In this case the vertex $(u^{''}, v-1)$ where $u^{''}\neq u, u^{'}$,
   is a common neighbor of $(u, v)$ and $(u^{'}, v^{'})$. Thus $d((u, v)$, $(u^{'}, v^{'}))=2$.  
  
  $iii)$ $u\neq u^{'}$ and $v\neq v^{'}$. Hence by $(i)$ and $(ii)$,  $d((u, v)$, $(u^{'}, v^{'}))=2$.
  
  Therefore $diam(\Gamma)=2$.
  
\begin{table}
\begin{center}
\caption{common neighbor between $(u, v), (u^{'}, v^{'})$ in $\Gamma=\Cay(\mathbb{Z}_{p}\times \mathbb{Z}_{p^{\alpha}q^{\beta}r^{\gamma}},\Phi)$}\label{tab:dia3}
\begin{tabular}{|c|c|c|}
\hline
$u=u^{'}, v\neq v^{'}$ & common neighbor & Comments \\
\hline
 $v, v^{'}\in \varphi_{p^{\alpha}q^{\beta}r^{\gamma}}$& $(u^{''}, pqr)$ &\\
 \hline 
 $v, v^{'}$ are multiples of $pqr$ & $(u^{''}, pqr-1)$ &\\
\hline
$v, v^{'}$ are multiples of $pq$ & $(u^{''}, r)$ & \\
\hline
$v, v^{'}$ are multiples of $pr$ & $(u^{''}, q)$ & \\
\hline
$v, v^{'}$ are multiples of $qr$ & $(u^{''}, p)$ & \\
\hline
$v, v^{'}$ are multiples of $p$ & $(u^{''}, qr)$ &\\
\hline
$v\neq v^{'}$ and each of them is & $(u^{''}, \dfrac{v+v^{'}}{2})$ & if $v-v^{'}\in \varphi_{p^{\alpha}q^{\beta}r^{\gamma}}$\\
 multiple of one of the prime factor & $(u^{''}, 2(v+v^{'}))$~~~~ & if $v-v^{'}\notin \varphi_{p^{\alpha}q^{\beta}r^{\gamma}}$\\
and both of them are even or odd &&\\
\hline
$v$ is multiple of $p$ & $(u^{''}, (\dfrac{v}{p})r+pqr)$ & if $\dfrac{v}{p}\in\varphi_{p^{\alpha}q^{\beta}r^{\gamma}}$\\
and $v^{'}$ is multiple of $pq$ & $(u^{''}, (\dfrac{v}{p})r+qr)$~~~~ & if ~$\dfrac{v}{p}$ be multiple of $p$\\
\hline
$v$ is multiple of $p$ & $(u^{''}, (\dfrac{v}{p})q+pqr)$ & if $\dfrac{v}{p}\in\varphi_{p^{\alpha}q^{\beta}r^{\gamma}}$\\
and $v^{'}$ is multiple of $pr$ & $(u^{''}, (\dfrac{v}{p})q+qr)$~~~~ & if ~$\dfrac{v}{p}$ be multiple of $p$\\
\hline
 $v$ is multiple of $p$& $(u^{''}, \dfrac{v+v^{'}}{2})$ & if $v,v^{'}$ be both even or odd\\
and $v^{'}$ is multiple of $qr$& $(u^{''}, 2(v+v^{'}))$~~~~ & if one of them be odd and other  \\
&& be even\\
\hline
$v$ is multiple of $p$ & $(u^{''}, \dfrac{v}{p}+pqr)$ & if ~$\dfrac{v}{p}$ be non-multiple of $p$\\
and $v^{'}$ is multiple of $pqr$ & $(u^{''}, \dfrac{v}{p}+qr)$~~~~ & if $\dfrac{v}{p}$ be multiple of $p$\\
\hline
$v$ is multiple of $p$ & $(u^{''}, (v+ v^{'})qr)$ & if $v,v^{'}$ be both even or odd\\
and $v^{'}\in\varphi_{p^{\alpha}q^{\beta}r^{\gamma}}$& $(u^{''}, v^{'}qr)$~~~~& if one of them be odd and other  \\
&& be even\\
\hline
$v, v^{'}$ are multiples of $q$ & $(u^{''}, pr)$ &\\
\hline
$v$ is multiple of $q$ & $(u^{''}, (\dfrac{v}{q})r+pqr)$ & if $\dfrac{v}{q}\in\varphi_{p^{\alpha}q^{\beta}r^{\gamma}}$\\
and $v^{'}$ is multiple of $pq$ & $(u^{''}, (\dfrac{v}{q})r+pr)$~~~~ & if ~$\dfrac{v}{q}$ be multiole of $q$\\
\hline
$v$ is multiple of $q$ & $(u^{''}, (\dfrac{v}{q})p+pqr)$ & if $\dfrac{v}{q}\in\varphi_{p^{\alpha}q^{\beta}r^{\gamma}}$\\
and $v^{'}$ is multiple of $qr$ & $(u^{''}, (\dfrac{v}{q})p+pr)$~~~~ & if ~$\dfrac{v}{q}$  be multiole of $q$\\
\hline
 $v$ is multiple of $q$& $(u^{''}, \dfrac{v+v^{'}}{2})$ & if $v,v^{'}$ be both even or odd\\
and $v^{'}$ is multiple of $pr$& $(u^{''}, 2(v+v^{'}))$~~~~ & if one of them be odd and other  \\
&& be even\\
\hline
$v$ is multiple of $q$ & $(u^{''}, \dfrac{v}{q}+pqr)$ & if ~$\dfrac{v}{q}$ be non-multiple of $q$\\
and $v^{'}$ is multiple of $pqr$ & $(u^{''}, \dfrac{v}{q}+pr)$~~~~ & if $\dfrac{v}{q}$ be multiple of $q$\\
\hline
$v$ is multiple of $q$ & $(u^{''}, (v+ v^{'})pr)$ & if $v,v^{'}$ be both even or odd\\
and $v^{'}\in\varphi_{p^{\alpha}q^{\beta}r^{\gamma}}$& $(u^{''}, v^{'}pr)$~~~~& if one of them be odd and other  \\
&& be even\\
\hline
\end{tabular}
\end{center}
\end{table}

\begin{table}
\begin{center}
\begin{tabular}{|c|c|c|}
\hline
$u=u^{'}, v\neq v^{'}$ & common neighbor & Comments \\
\hline
$v, v^{'}$ are multiple of $r$ & $(u^{''}, pq)$ & \\
\hline
$v$ is multiple of $r$ & $(u^{''}, (\dfrac{v}{r})q+pqr)$ & if $\dfrac{v}{r}\in\varphi_{p^{\alpha}q^{\beta}r^{\gamma}}$\\
and $v^{'}$ is multiple of $pr$ & $(u^{''}, (\dfrac{v}{r})q+pq)$~~~~ & if ~$\dfrac{v}{r}$ be multiole of $r$\\
\hline
$v$ is multiple of $r$ & $(u^{''}, (\dfrac{v}{r})p+pqr)$ & if $\dfrac{v}{r}\in\varphi_{p^{\alpha}q^{\beta}r^{\gamma}}$\\
and $v^{'}$ is multiple of $qr$ & $(u^{''}, (\dfrac{v}{r})p+pq)$~~~~ & if ~$\dfrac{v}{r}$  be multiole of $r$\\
\hline
 $v$ is multiple of $r$& $(u^{''}, \dfrac{v+v^{'}}{2})$ & if $v,v^{'}$ be both even or odd\\
and $v^{'}$ is multiple of $pq$& $(u^{''}, 2(v+v^{'}))$~~~~ & if one of them be odd \\
\hline
$v$ is multiple of $r$ & $(u^{''}, \dfrac{v}{r}+pqr)$ & if ~$\dfrac{v}{r}$ be non-multiple of $r$\\
and $v^{'}$ is multiple of $pqr$ & $(u^{''}, \dfrac{v}{r}+pq)$~~~~ & if $\dfrac{v}{r}$ be multiple of $r$\\
\hline
$v$ is multiple of $r$ & $(u^{''}, (v+ v^{'})pq)$ & if $v,v^{'}$ be both even or odd\\
and $v^{'}\in\varphi_{p^{\alpha}q^{\beta}r^{\gamma}}$& $(u^{''}, v^{'}pq)$~~~~& if one of them be odd  \\
\hline
$v$ is multiple of $pq$ & $(u^{''}, \dfrac{v+v^{'}}{p})$ & if ~$\dfrac{v+v^{'}}{p}\in\varphi_{p^{\alpha}q^{\beta}r^{\gamma}}$ \\
and $v^{'}$ is multiple of $pr$ & $(u^{''}, \dfrac{v+v^{'}}{p}+qr)$~~~~ & if $\dfrac{v+v^{'}}{p}$ be multiple of $p$\\
\hline
$v$ is multiple of $pq$ & $(u^{''}, \dfrac{v+v^{'}}{q})$ & if ~$\dfrac{v+v^{'}}{q}\in \varphi_{p^{\alpha}q^{\beta}r^{\gamma}}$\\
and $v^{'}$ is multiple of $rq$ & $(u^{''}, \dfrac{v+v^{'}}{q}+pr)$~~~~ & if $\dfrac{v+v^{'}}{q}$ be multiple of $q$\\
\hline
$v$ is multiple of $pr$ & $(u^{''}, \dfrac{v+v^{'}}{r})$ & if ~$\dfrac{v+v^{'}}{r}\in \varphi_{p^{\alpha}q^{\beta}r^{\gamma}}$\\
and $v^{'}$ is multiple of $rq$ & $(u^{''}, \dfrac{v+v^{'}}{r}+pq)$~~~~ & if $\dfrac{v+v^{'}}{r}$ be multiple of $r$\\
\hline
$v=kpq,k\in\mathbb{Z},v^{'}\in\varphi_{p^{\alpha}q^{\beta}r^{\gamma}}$ and& $(u^{''}, \dfrac{v+v^{'}}{2})$ & if ~$v-v^{'}\in\varphi_{p^{\alpha}q^{\beta}r^{\gamma}}$ \\ $v,v^{'}$ are both even or odd & $(u^{''}, v-v^{'})$~~~~ & if $v-v^{'}\notin\varphi_{p^{\alpha}q^{\beta}r^{\gamma}}$\\
one of the $v$ or $v^{'}$ is odd & $(u^{''}, v^{'}r)$&if $v-v^{'}\notin\varphi_{p^{\alpha}q^{\beta}r^{\gamma}}$\\
\hline
$v=kpr,k\in\mathbb{Z},v^{'}\in\varphi_{p^{\alpha}q^{\beta}r^{\gamma}}$and & $(u^{''}, \dfrac{v+v^{'}}{2})$ & if ~$v-v^{'}\in\varphi_{p^{\alpha}q^{\beta}r^{\gamma}}$ \\
$v,v^{'}$ are both even or odd & $(u^{''}, v-v^{'})$~~~~ & if $v-v^{'}\notin\varphi_{p^{\alpha}q^{\beta}r^{\gamma}}$\\
one of the $v$ or $v^{'}$ is odd & $(u^{''}, v^{'}q)$&if $v-v^{'}\notin\varphi_{p^{\alpha}q^{\beta}r^{\gamma}}$\\
\hline
$v=kqr, k\in\mathbb{Z},v^{'}\in\varphi_{p^{\alpha}q^{\beta}r^{\gamma}}$ and& $(u^{''}, \dfrac{v+v^{'}}{2})$ & if ~$v-v^{'}\in\varphi_{p^{\alpha}q^{\beta}r^{\gamma}}$ \\
$v,v^{'}$ are both even or odd & $(u^{''}, v-v^{'})$~~~~ & if $v-v^{'}\notin\varphi_{p^{\alpha}q^{\beta}r^{\gamma}}$\\
one of the $v$ or $v^{'}$ is odd & $(u^{''}, v^{'}p)$&if $v-v^{'}\notin\varphi_{p^{\alpha}q^{\beta}r^{\gamma}}$\\
\hline
 &  & by Proposition 4.5\cite{me}, $\lambda=4$  \\ 
$v=kpqr, k\in\mathbb{Z}$ and $v^{'}$& $(u^{''}, +4)$ or $(u^{''}, -4)$ & then $pqr$ is adjacent by $\pm 4$ \\
  is multiple of $pq$ or $pr$ or $qr$& & and $pq, pr, qr$ are\\ 
& & adjacent by $+4$ or $-4$\\
\hline
$v$ is multiple of $pqr$ & $(u^{''}, \dfrac{v+v^{'}}{2})$ & if $v,v^{'}$ be both even or odd\\
and $v^{'}\in \varphi_{p^{\alpha}q^{\beta}r^{\gamma}}$  & $(u^{''}, 2(v+v^{'}))$~~~~ & if one of them be odd\\
\hline
\end{tabular}
\end{center}
\end{table}
\end{proof}
\end{prop}
\begin{prop}\label{t6}
Let $\Gamma=\Cay(\mathbb{Z}_{p}\times\mathbb{Z}_{p^{\alpha}q^{\beta}r^{\gamma}}, \Phi)$,
where $p, q, r \geq 3 $ and $\alpha, \beta, \gamma\geq 1$.
Then $\gamma_{t}(\Gamma)$ and $\gamma_{c}(\Gamma)$ is given by Table~\ref{tab:4}.
\begin{table}
\begin{center}
\caption{$\gamma_{t}(\Cay(\mathbb{Z}_{p}\times \mathbb{Z}_{p^{\alpha}q^{\beta}r^{\gamma}}, \Phi))$ and $\gamma_{c}(\Cay(\mathbb{Z}_{p}\times \mathbb{Z}_{p^{\alpha}q^{\beta}r^{\gamma}}, \Phi))$ where $p, q, r \geq 3 $}
\label{tab:4}
\begin{tabular}{|c|c|c|}
\hline
$\Gamma$& $\gamma_{t}(\Gamma),\gamma_{c}(\Gamma) $ & Comments \\
\hline
$\Cay(\mathbb{Z}_{p}\times \mathbb{Z}_{p^{\alpha}q^{\beta}r^{\gamma}}, \Phi)$ & $6\leq\gamma_{t}(\Gamma), \gamma_{c}(\Gamma) \leq 8$
 & one of the prime factors is $3$\\
\hline
$\Cay(\mathbb{Z}_{p}\times \mathbb{Z}_{p^{\alpha}q^{\beta}r^{\gamma}}, \Phi)$ & $5$ & $p, q, r \geq 5 $  \\
\hline
\end{tabular}
\end{center}
\end{table}
\begin{proof}
By using the \cite[Proposition 4.5]{me}, $D, D^{'}, D^{''}, D^{'''}$ are minimal dominating sets for various cases in this graph.
 Clearly the subgraphs generated by $D$, $D^{'}$, $D^{''}$ and $D^{'''}$ are all connected. Therefore $\gamma(\Gamma)=\gamma_{t}(\Gamma)=\gamma_{c}(\Gamma).$
\end{proof}
\end{prop}

As an immediate consequence of  Lemma~\ref{t4} and Propositions~\ref{t5},~\ref{t6},
we have the following theorem.

\begin{thm}
Let $\Gamma=\Cay(\mathbb{Z}_{p}\times\mathbb{Z}_{p^{\alpha}q^{\beta}r^{\gamma}}, \Phi)$, where
 $p, q, r \geq 2$ and $\alpha, \beta, \gamma \geq 1$. Then
 $\gamma_{t}(\Gamma)$and $\gamma_{c}(\Gamma)$ is given by Table \ref{tab:5}.

\begin{table}
\begin{center}
\caption{$\gamma_{t}(\Cay(\mathbb{Z}_{p}\times\mathbb{Z}_{p^{\alpha}q^{\beta}r^{\gamma}}, \Phi))$ and $\gamma_{c}(\Cay(\mathbb{Z}_{p}\times\mathbb{Z}_{p^{\alpha}q^{\beta}r^{\gamma}}, \Phi))$ where $\alpha, \beta, \gamma\geq 1$}\label{tab:5}
\begin{tabular}{|c|c|c|c|}
\hline
$\Gamma$& $\gamma_{t}(\Gamma)$ & $\gamma_{c}(\Gamma)$ & Comments \\
\hline
$ \Cay(\mathbb{Z}_{2}\times \mathbb{Z}_{2^{\alpha}q^{\beta}r^{\gamma}}, \Phi)$ & $12$ & does not exist & \\
\hline
$\Cay(\mathbb{Z}_{p}\times \mathbb{Z}_{2^{\alpha}p^{\beta}r^{\gamma}}, \Phi)$ & $10$ & $12$ & one of the prime factors is $3$\\
&&& $p=3$\\
\hline
$\Cay(\mathbb{Z}_{p}\times \mathbb{Z}_{2^{\alpha}p^{\beta}r^{\gamma}}, \Phi)$ & $10$ & $10$ & one of the prime factors is $3$\\
&&& $p\geq 5 $  \\
\hline
$ \Cay(\mathbb{Z}_{p}\times \mathbb{Z}_{2^{\alpha}p^{\beta}r^{\gamma}}, \Phi)$ & $8$ & $8$ & $p, r\geq 5$\\
\hline
$ \Cay(\mathbb{Z}_{p}\times \mathbb{Z}_{p^{\alpha}q^{\beta}r^{\gamma}}, \Phi)$ & $6\leq\gamma_{t}(\Gamma)\leq 8$ &  $6\leq\gamma_{c}(\Gamma)\leq 8$ & one of the prime factors is $3$\\
\hline
$\Cay(\mathbb{Z}_{p}\times \mathbb{Z}_{p^{\alpha}q^{\beta}r^{\gamma}}, \Phi)$ & $5$ & $5$ & $p, q, r \geq 5 $\\
\hline
\end{tabular}
\end{center}
\end{table}
\end{thm}
As an immediate consequence of  Lemmas~\ref{d2},~\ref{tdia} and Propositions~\ref{d3},~\ref{d4},~\ref{d6},~\ref{d7},
we have the following theorem.
\begin{thm}
Let $\Gamma=\Cay(\mathbb{Z}_{p}\times\mathbb{Z}_{m}, \Phi)$, where $m\in\{p^{\alpha}q^{\beta}, p^{\alpha}q^{\beta}r^{\gamma}\}$. Then
\item[1)] $diam(\Gamma)=3$ where one of the prime factors be $2$.
\item[1)] $diam(\Gamma)=2$ where $p, q, r\geq 3$.
\end{thm}
\begin{remark}
Let $ p_{1}, p_{2}, \ldots, p_{k}$ be consecutive prime numbers, $p_{1}=3$, $\alpha, \alpha_{1}, \alpha_{2}, \ldots, \alpha_{k}$ are
positive integers, $\alpha\geq 2$, and $\Phi=\varphi_{2}\times \varphi_{2^{\alpha}p_{1}^{\alpha_{1}}p_{2}^{\alpha_{2}}
 \ldots p_{k}^{\alpha_{k}}}$. Then by \cite[Theorem 4.7]{me}, we have $\gamma(\Cay(\mathbb{Z}_{2}\times \mathbb{Z}_{2^{\alpha}p_{1}^{\alpha_{1}}p_{2}^{\alpha_{2}}
 \ldots p_{k}^{\alpha_{k}}}, \Phi))\geq 4k+4$. Therefore 
 $\gamma_{t}(\Cay(\mathbb{Z}_{2}\times \mathbb{Z}_{2^{\alpha}p_{1}^{\alpha_{1}}p_{2}^{\alpha_{2}}
 \ldots p_{k}^{\alpha_{k}}}, \Phi))\geq 4k+4$. Since $\Gamma$ is a disconnected graph the {\it $\gamma_{c}$-set} does not exist for $\Gamma$. 

\end{remark}

%\begin{acknowledgements}
%\end{acknowledgements}



\end{document}